\documentclass[paper,twocolumn,twoside]{geophysics}

\usepackage{amsmath}

\begin{document}

\title{A Novel Immersed Boundary Approach for Irregular \\ 
       Topography with Acoustic Wave Equations}

\renewcommand{\thefootnote}{\fnsymbol{footnote}} 

\address{
\footnotemark[1]Imperial College London, Department of Earth Science and Engineering, London, UK
\footnotemark[2]Devito Codes, London, UK}
\author{Edward Caunt\footnotemark[1], Rhodri Nelson\footnotemark[1], Fabio Luporini\footnotemark[2], Gerard Gorman\footnotemark[1]}

\footer{Example}
\lefthead{Caunt et al.}
\righthead{Immersed boundary topography approach}


\begin{abstract}
  Irregular terrain has a pronounced effect on the propagation of seismic and acoustic wavefields but is not straightforwardly reconciled with structured finite-difference (FD) methods used to model such phenomena. Methods currently detailed in the literature are generally limited in scope application-wise or non-trivial to apply to real-world geometries. With this in mind, a general immersed boundary treatment capable of imposing a range of boundary conditions in a relatively equation-agnostic manner has been developed, alongside a framework implementing this approach, intending to complement emerging code-generation paradigms. The approach is distinguished by the use of N-dimensional Taylor-series extrapolants constrained by boundary conditions imposed at some suitably-distributed set of surface points. The extrapolation process is encapsulated in modified derivative stencils applied in the vicinity of the boundary, utilizing hyperspherical support regions. This method ensures boundary representation is consistent with the FD discretization: both must be considered in tandem. Furthermore, high-dimensional and vector boundary conditions can be applied without approximation prior to discretization. A consistent methodology can thus be applied across free and rigid surfaces with the first and second-order acoustic wave equation formulations. Application to both equations is demonstrated, and numerical examples based on analytic and real-world topography implementing free and rigid surfaces in 2D and 3D are presented.
\end{abstract}

\section{Introduction}

Irregular topography introduces considerable complexity to geophysical wavefield models, creating complex path effects which turn clean, defined reflections into cascades of overlapping arrivals. Interactions with topography cause waves to diffract and scatter, focus and defocus \citep{Takemura2015, Reinoso1996, Boore1972, Griffiths1979}: effects that must be encapsulated by any attempt to simulate the behaviour of such wavefields \citep{Borisov2018}. Whilst early numerical experiments demonstrated its capacity to markedly affect recorded data \citep{Boore1973}, the proliferation of wave-equation-based workflows, most notably full-waveform inversion (FWI) and reverse-time migration (RTM), has placed a sharpened focus on understanding topographic effects.

Topographic effects have been explored in contexts including FWI \citep{Bleibinhaus2009}, understanding seismic wave scattering \citep{Takemura2015}, infrasound propagation problems \citep{Kim2014, Fee2021}, and is crucial for the emerging field of teleseismic FWI \citep{Monteiller2013, Monteiller2015}. It is recognised that poor topography implementation adversely affects model accuracy (e.g. \citealp{Monteiller2013, Li2020, Borisov2018}), often severely \citep{Nuber2016} although this may be acceptable in some applications \citep{Bleibinhaus2009}. When applied to imaging problems, unrealistic propagation paths result in artefacts in the processed image \citep{Bleibinhaus2009, Nuber2016, Borisov2018} which risk being mistaken for real geological features, whilst models of ground motion have been found to significantly underestimate local amplification factors if topography is omitted \citep{Reinoso1996}.

Unstructured meshes suitably conformed to the topography offer an immediate solution. Finite-element methods (FEMs) have been demonstrated as effective for modelling wave propagation in the presence of variable topography (e.g. \citealp{Zhebel2014, Dupros2010, Borisov2018, Liu2014, Mulder2016}), and have been applied to include complex terrain in models of earthquake peak ground acceleration \citep{Galis2008}, in shallow structure characterization applications \citep{Romdhane2011}, and to imaging problems such as full waveform inversion (FWI) \citep{Roberts2021, Shin2013, Monteiller2015}. Unstructured finite-difference methods (FDMs) have also been developed \citep{Takekawa2015, Martin2015}, and successfully used to perform FWI on synthetic datasets \citep{Wu2021}.

However, structured FD methods dominate in geophysical wavefield modelling and processing applications, including seismic reverse-time migration (RTM) \citep{Fletcher2009}, full-waveform inversion (FWI) \citep{Warner2013}, and source localisation \citep{Mckee2014}. This is not without reason: they eschew potentially cumbersome grid-generation algorithms \citep{Brehm2013}, which become increasingly problematic at large scales \citep{Slotnick2014} and require apriori knowledge of seismic velocities \citep{Roberts2021a, Roberts2021}. The requisite information is potentially unavailable in practice or is iteratively updated in outer-loop problems, necessitating repeated mesh adaptation. Furthermore, the geometry of geological structures is prone to generating ill-conditioned sliver elements as units get pinched out resulting in an unreasonably small timestep \cite{Roberts2021}. Automated mesh generation for seismic wave propagation in models containing arbitrary horizons remains an open problem, and thus unstructured approaches are rarely used in production. Structured FD methods are simple to implement, with relatively-low computational footprints \citep{Liu2009}, even before considering the suite of known optimizations \citep{Luporini2019, Louboutin2019}. However, in the presence of irregular topography, accurately representing the sharp, uneven material discontinuity on regular grids can be problematic \citep{Zeng2012, Mulder2017, Gao2015}. Ideally, one would want to accurately represent complex, curvilinear topography as a sharp interface whilst retaining the advantages of structured grids.

Air- or vacuum-layer approaches, achieved through a low-impedence layer at the surface \citep{Schultz1997, Zeng2012} trade accuracy for simplicity in mimicking the free surface. Satisfactory results can be obtained with careful implementation \citep{Zeng2012}, although error analyses and numerical experiments indicate approximate sub-second-order convergence in space, regardless of interior discretization \citep{Zhebel2014, Symes2009} and potentially egregious \citep{Graves1996, Zahradnik1993} error when applied to complex geometries. Suppression of spurious scattering requires heavy oversampling \citep{Bohlen2006}, and smoothing or stabilisation routines are often required to stabilise the surface treatment \citep{Bartel2000, Zeng2012, Vieira2018}. Whilst such approaches can be improved with locally-refined subgrids \citep{Oprsal1999, Tavelli2019}, such refinement is non-trivial \citep{Lai2000, Goldstein1993}, often requiring careful filtering and interpolation routines \citep{Zhang2013a}. Staircased image methods marginally improve on this \citep{Robertsson1996}: variations have been explored by several authors (e.g. \citet{Robertsson1996, Boore1972, Hayashi2001, Ohminato1997, Ripperger2003}), used to study macro-scale topographic scattering effects \citep{Takemura2015, Nakamura2012} and volcanic source location \citep{Kim2014, Fee2021}. The stepped boundary still generates diffractions \citep{Muir1992, Hu2016}, and a time shift proportional to the difference between the interface and computational grids \citep{Symes2009}. A more accurate image method can be achieved via coordinate transform (e.g. \citealp{Petersson2015}) such that the surface is a horizontal plane in the iteration space \citep{Zhang2006, Hestholm2002}. This approach is widely adaptable, and various wave equations have been solved with such schemes (e.g. \citealp{Zhang2006, Zhang2012, Hestholm1994, Hestholm1999, Hestholm2000, Sun2017, Puente2014}), demonstrating a high degree of accuracy, on par with FEMs \citep{Zhang2012}. However, a smooth conformal mapping may be challenging to obtain or yield locally small cells, limiting the maximum timestep \citep{Shragge2014} in much the same manner as the aforementioned sliver elements in FEM.

An alternative approach, popular in computational fluid dynamics (CFD) contexts (see \citealp{Mittal2008, Seo2011, Dong2010, Vargas2008} for examples), is the immersed-boundary method. This approach embeds a curvilinear interface within a Cartesian grid by locally modifying the finite-difference operators in the vicinity of the boundary \citep{Brehm2013, Brehm2015, Mulder2017}. This approach has seen several variations for second-order acoustic wave equation formulations \citep{Zhang2013, Mulder2017, Li2020}, with some extending to first-order formulations \citep{Mulder2017a, Hu2016}, and even the isotropic elastic wave equation with some success \citep{Lombard2008, Almuhaidib2015, Gao2015}. The standard Cartesian grid and equations are retained and complex geometries including sharp edges and concavities can be represented.

To date, most applications of such schemes in geophysical wave modelling scenarios have made use of problem-specific approximations prior to discretization to impose boundary conditions. For example, imposing a suitable 1D approximation of the boundary conditions. It has been repeatedly demonstrated that such approximations yield satisfactory results \citep{Hu2016, Mulder2017, Almuhaidib2015, Li2022} and their usage is widespread. Whilst such approaches offer intuitive parallels with conventional image methods, the condition imposed may not, even at its limit, strictly equal the true boundary condition, instead merely being suitably similar. This is particularly apparent for vector boundary conditions such as those applied to the velocity fields at the free surface. \citet{Hu2016} addressed this by decomposing velocity field components into tangential and radial components in a local cylindrical coordinate system, whilst \citet{Mulder2017a} assumed locally vertical or horizontal boundaries, incurring numerical error.

Such approximations can be suitable on a case-by-case basis, but are reliant upon domain knowledge and do not necessarily generalise straightforwardly: it is not immediately clear as to their applicability across multiple wave equations. Indeed, these approximations may be inherently limiting, precluding adaption to more complex physics (wave equations featuring transverse isotropy or elastic wave equations for example). Given the trajectory towards increasingly complex physics and geometries, this paper sets out to develop a generalised immersed boundary approach, such that a consistent methodology can be applied in a relatively problem-agnostic manner.

This approach aligns with current trends in geophysical modelling: domain-specific languages (DSLs) and automatic code generation are increasingly prevalent, with projects such as Devito \citep{Luporini2019, Louboutin2019} and Firedrake \citep{Rathgeber2016} leveraging high-level abstractions to generate low-level FDM and FEM solver kernels from symbolic partial differential equations (PDEs). Abstracting low-level aspects of implementation achieves a separation of concerns between the overarching problem and its underlying implementation \citep{Rathgeber2016}, enabling domain specialists to focus on the application level. High-level interfaces substantially reduce development time \citep{Louboutin2019}, whilst sophisticated optimization routines in the lowering process yield high-performance, portable code \citep{Luporini2019, Louboutin2019}. Such approaches hinge on a general method for solving a given problem so that suitable abstractions can be developed; this cannot be easily achieved with the application-specific immersed boundary implementations developed to date.

Developing a generalised mathematical approach enables abstraction, which in turn lends itself to automation and thus code generation. Furthermore, direct generalisation to other wave equations widens potential application, spanning a range of geophysical problems. The generalised nature of the method presented greatly simplifies inclusion of immersed boundaries in geophysical models, as a suitable treatment can be devised according to the prescribed formula to suit the physical problem at hand.

This paper is structured as follows: an immersed boundary approach supporting the imposition of multidimensional and vector boundary conditions is outlined, followed by its application to a selection of equations and boundary conditions of interest. This mathematical approach forms the backbone of Schism: a plugin for Devito used to implement the examples shown in this paper. The numerical accuracy of the approach proposed is explored through convergence testing of the resultant treatment, and several geophysically-relevant test cases based on the first and second-order formulations of the acoustic wave equation are explored. Note that whilst the demonstrations within this paper focus solely on variants of the acoustic wave equation, the method is nominally equally suited to pseudoacoustic wave equations featuring vertical and tilted transverse isotropy (VTI and TTI), alongside elastic wave equations. Further exploration of these areas is planned in subsequent publications.

\section*{General Approach for Constructing Immersed Boundary Derivative Operators}
\renewcommand{\figdir}{Theory} 
\subsection{Constructing Finite-Difference Approximations}
In this section, the discretization of continuous fields in the vicinity of non-grid-conforming boundaries is outlined. The reduction of continuous functions to a discrete set of values is key to all manner of numerical models. In the case of FDMs, these points constitute nodes on a grid, and updating these discrete values according to some governing equation approximates the evolution of the underlying continuous function. Function values located at these points constrain a basis with some specified error relative to the continuous function, diminishing as resolution increases.

Understanding the evolution of these functions necessitates the calculation of derivatives; to this end the aforementioned basis is used to form approximations of the derivative operators, discretizing the continuous equation. As with the basis, these will have some error relative to their true value. To outline how these approximations are made in the typical case, consider some function $f$ in a 1D space; approximating with an $M$th-order Taylor series expansion at some point $x_0$ yields
\begin{equation}
    f(x) = \sum_{m=0}^M \frac{(x-x_0)^{m}}{m!} \frac{\partial^{m}f(x_0)}{\partial x^{m}} + O((x-x_0)^{M+1}) \mathrm{,}
\end{equation}
concisely represented as
\begin{equation}
    \mathbf{a}\cdot\boldsymbol{\delta} = f(\mathbf{\mathrm{x}}) \mathrm{,}
    \label{eq:accmaths:dot_product}
\end{equation}
where
\begin{equation}
    \mathbf{a}^T = \left(1, (x-x_0), ..., \frac{(x-x_0)^M}{M!}\right)\mathrm{,}
\end{equation}
and
\begin{equation}
    \boldsymbol{\delta}^T = \left(f(x_0), \frac{\partial f}{\partial x}(x_0), ..., \frac{\partial^M f}{\partial x^M}(x_0)\right)\mathrm{.}
\end{equation}
For some even $M$, expansions at $M+1$ discrete points, labelled $x_{-M/2}$ through $x_{M/2}$, enables the formation of the linear system
\begin{equation}
    \mathbf{A}\boldsymbol{\delta} =
    \begin{pmatrix}
        f(x_{-M/2}) \\
        \vdots \\
        f(x_{M/2})
    \end{pmatrix} \mathrm{.}
    \label{eq:accmaths:standard_linear_system}
\end{equation}
To illustrate, for $M=2$, the system is as follows
\begin{equation}
    \begin{pmatrix}
    1 & x_{-1}-x_0 & \frac{(x_{-1}-x_0)^2}{2} \\
    1 & 0 & 0 \\
    1 & x_1-x_0 & \frac{(x_1-x_0)^2}{2}
    \end{pmatrix}
    \begin{pmatrix}
        f(x_0) \\
        \frac{\partial f}{\partial x}(x_0) \\
        \frac{\partial^2 f}{\partial x^2}(x_0)      
    \end{pmatrix}
    =
    \begin{pmatrix}
        f(x_{-1}) \\
        f(x_0) \\
        f(x_1)    
    \end{pmatrix}\mathrm{,}
\end{equation}
which can be rearranged in the form
\begin{equation}
    \mathbf{A}^{-1}
    \begin{pmatrix}
        f(x_{-M/2}) \\
        \vdots \\
        f(x_{M/2})
    \end{pmatrix}
    = \boldsymbol{\delta} \mathrm{,}
\end{equation}
thereby obtaining derivatives as some weighted sum of discrete function values: a stencil. With the assumption of a regular grid, relative positions of points become fixed regardless of where the derivative is being taken, yielding
\begin{equation}
    \begin{pmatrix}
        1 & -1 & \frac{1}{2} \\
        1 & 0 & 0 \\
        1 & 1 & \frac{1}{2}
    \end{pmatrix}
    \begin{pmatrix}
        f(x_0) \\
        \Delta x\frac{\partial f}{\partial x}(x_0) \\
        \Delta x^2\frac{\partial^2 f}{\partial x^2}(x_0)      
    \end{pmatrix}
    =
    \begin{pmatrix}
        f(x_{-1}) \\
        f(x_0) \\
        f(x_1)   
    \end{pmatrix}\mathrm{,}
\end{equation}
for the above case, the inverse being
\begin{equation}
    \begin{pmatrix}
        0 & 1 & 0 \\
        -\frac{1}{2} & 0 & \frac{1}{2} \\
        1 & -2 & 1
    \end{pmatrix}
    \begin{pmatrix}
        f(x_{-1}) \\
        f(x_0) \\
        f(x_1)     
    \end{pmatrix}
    =
    \begin{pmatrix}
        f(x_0) \\
        \Delta x\frac{\partial f}{\partial x}(x_0) \\
        \Delta x^2\frac{\partial^2 f}{\partial x^2}(x_0)       
    \end{pmatrix} \mathrm{.}
\end{equation}
One may observe that the leftmost matrix contains weights for FD stencils of every derivative order up to that of the basis.

For higher dimensions, a suitably higher-dimensional polynomial is required, formed as a product of per-dimension 1D series of the form
\begin{equation}
    f(x) = \sum_{m=0}a_m x^m \mathrm{,}
\end{equation}
in 3D yielding
\begin{equation}
    f(\mathbf{\mathrm{x}}) = \sum_{m=0}\sum_{n=0}\sum_{l=0} a_{mnl} x^m y^n z^l \mathrm{.}
\end{equation}
Given this, the process outlined above is extensible to higher dimensions as in \citet{Takekawa2015}. In the pursuit of consistency and accuracy, it is crucial to maintain consistent error throughout any numerical scheme; thus polynomials used in such schemes must be of an order matching the spatial discretization employed. To this end, for an order $M$ spatial discretization, the 3D expansion is truncated as follows:
\begin{equation}
    f(\mathbf{\mathrm{x}}) = \sum_{m=0}^{M_x}\sum_{n=0}^{M_y}\sum_{l=0}^{M_z} a_{mnl} x^m y^n z^l + O(|\Delta\mathbf{\mathrm{x}}|^{(M+1)}) \mathrm{,}
\end{equation}
subject to
\begin{equation}
    M_x + M_y + M_z = M \mathrm{,}
    \label{eq:accmaths:truncation}
\end{equation}
thus removing any terms with order greater than $M$.

An N-dimensional Taylor series is a robust candidate for such applications, retaining similarities to the aforementioned 1D case and being suitably problem-agnostic. In 3D, this is given by
\begin{equation}
    f(\mathbf{\mathrm{x}}) = \sum_{m=0}\sum_{n=0}\sum_{l=0} \frac{(x-x_0)^m(y-y_0)^n(z-z_0)^l}{m!n!l!} \delta_{m,n,l})\mathrm{,}
    \label{eq:accmaths:3D_taylor}
\end{equation}
where
\begin{equation}
    \delta_{m,n,l} = \frac{\partial^{m+n+l}}{\partial x^m \partial y^n \partial z^l}f(\mathbf{\mathrm{x}}_0) \mathrm{.}
\end{equation}
Truncating this expansion as in Equation \ref{eq:accmaths:truncation}, the polynomial expansion at points within the support region can be represented as a matrix-vector multiplication, as in Equation \ref{eq:accmaths:dot_product}. Reflecting the higher dimensionality and corresponding increased number of terms, $\mathbf{a}$ is given by
\begin{equation}
\begin{aligned}
    \mathbf{a}^T = \bigg( & 1, (x-x_0), ..., (z-z_0), \frac{(x-x_0)^2}{2}, \\
    & (x-x_0)(y-y_0), ..., \frac{(z-z_0)^2}{2}, ..., \\
    & \frac{(x-x_0)^M}{M!}, ..., \frac{(z-z_0)^M}{M!}\bigg) \mathrm{,}
\end{aligned}
\end{equation}
whilst
\begin{equation}
\begin{aligned}
    \boldsymbol{\delta}^T = \bigg( & 1, \frac{\partial}{\partial x},..., \frac{\partial}{\partial z}, \frac{\partial^2}{\partial x^2}, \\
    & \frac{\partial^2}{\partial x \partial y}, ..., \frac{\partial^2}{\partial z^2}, ..., \\
    & \frac{\partial^M}{\partial x^M},..., \frac{\partial^M}{\partial z^M} \bigg) f(\mathbf{\mathrm{x}}_0)\mathrm{.}
\end{aligned}
\end{equation}

As in the 1D case, expansions taken at some set of points are used to form a linear system. However, with increased dimensionality comes increased flexibility regarding the distribution of these points. Consider function values discretized at a set of points distributed in an arbitrary manner. In one dimension it is reasonably simple to construct a polynomial expansion of any desired order to fit this data, provided the number of linearly-independent data points equals or exceeds the coefficients in the truncated expansion. However, when constraining the higher-dimensional basis, this linear independence requires suitable point distribution.

\plot{circular_support.png}{width=\columnwidth}
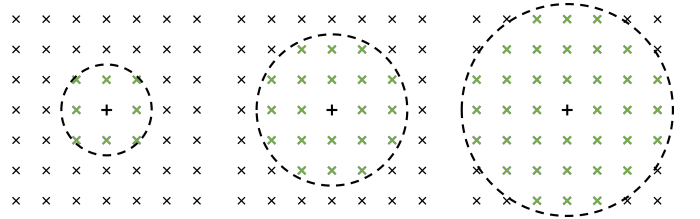
{A selection of 2D support region footprints with radii of 1.5, 2.5, and 3.5 respectively. The vertical black cross designates the stencil centre (the position at which the stencil is applied), green crosses are interior points, and the dotted black line shows the extent of the support region.
}

A hyperspherical support region footprint, similar to that used by \citet{Takekawa2015}, is proposed as a suitable choice for a large range of boundary conditions. Extensible to any number of dimensions, this choice of support region has a defined centre and can be straightforwardly expanded and contracted as necessary without introducing drastic behavioural changes. Furthermore, any given number of points within the support region will be located as close to the stencil centre as possible. Whilst this work limits itself to Cartesian structured grids, such support regions are nominally equally suited to structured and unstructured data as no assumption is required regarding the underlying data structure. In this case, points can be located by index, removing the need for nearest-neighbour searches to find points within the support region. Support regions of various radii are shown in Figure \ref{fig:circular_support.png}: a radius of $(M+1)/2 \Delta x$ ensures sufficient information to constrain a polynomial basis on an even, structured grid.

\plot{matrix_free_space_support.png}{width=\columnwidth}
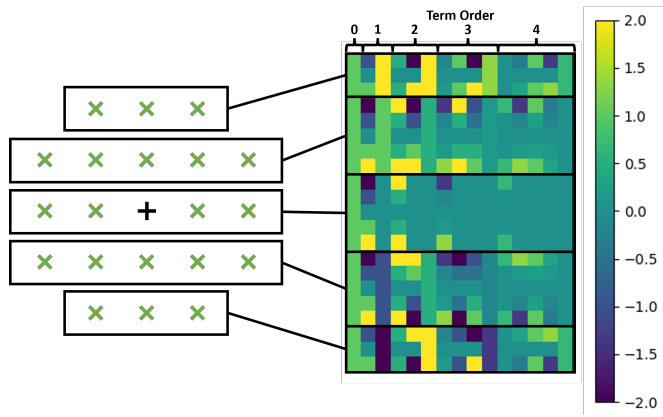
{A 2D support region footprint for a 4th-order basis consisting of points $\mathbf{X}$ and its correspondence to the matrix $\mathbf{A}$ (the structure of which is shown on the right). The stencil points in each box correspond to the matrix rows in the connected box, containing coefficients of the polynomial expansion at that point. The order of derivatives to which coefficients correspond increase from left to right. Colours within the matrix plot correspond to the values of those elements, with yellows and blues representing positive and negative values respectively. The vertical black cross designates the stencil centre: the position at which the stencil is applied and expansion point for the polynomial.
}
Whilst the approach taken readily extends to any number of dimensions, the proceeding example considers the 2D case for clarity and brevity. Correspondence between a 2D support region for a 4th-order discretization and $\mathbf{A}$ is shown in Figure \ref{fig:matrix_free_space_support.png}. Note that this matrix is tall, and thus the linear system
\begin{equation}
    \mathbf{A}\boldsymbol{\delta} = f(\mathbf{X})\mathrm{,}
\end{equation}
where
\begin{equation}
    \begin{aligned}
    \boldsymbol{\delta}^T = \bigg( & 1, \Delta x \frac{\partial}{\partial x}, \Delta y\frac{\partial}{\partial y}, \Delta x^2 \frac{\partial^2}{\partial x^2}, \Delta x \Delta y \frac{\partial^2}{\partial x \partial y}, \Delta y^2 \frac{\partial^2}{\partial y^2}, \\
    & \Delta x^3 \frac{\partial^3}{\partial x^3}, \Delta x^2 \Delta y \frac{\partial^3}{\partial x^2 \partial y}, \Delta x \Delta y^2 \frac{\partial^3}{\partial x \partial y^2}, \Delta y^3 \frac{\partial^3}{\partial y^3}, \\
    & \Delta x^4 \frac{\partial^4}{\partial x^4}, \Delta x^3 \Delta y\frac{\partial^4}{\partial x^3 \partial y}, \Delta x^2 \Delta y^2 \frac{\partial^4}{\partial x^2 \partial y^2}, \\
    & \Delta x \Delta y^3 \frac{\partial^4}{\partial x \partial y^3}, \Delta y^4 \frac{\partial^4}{\partial y^4} \bigg) f(\mathbf{\mathrm{x}}_0)\mathrm{.}
\end{aligned}
\end{equation}
will be overdetermined. This system is solved with the Moore-Penrose pseudoinverse to obtain the derivative vector. Inverting this system yields expressions for FD stencils using the aforementioned 2D basis.

\subsection{Imposition of Boundaries}
Consider some continuous scalar function within a domain. This could be a standalone scalar function or individual component of some vector field. In the context of wavefield modelling, this may represent pressure, a component of the particle velocity vector, or even some stress tensor component. The evolution of a field is dictated not only by the governing system of partial-differential-equations, but boundary conditions imposed. In the case of seismic wavefield modelling, these demark a half-space, bounded at the top by the ground-air interface. By imposing suitable conditions on this interface, the wavefield can be made to interact with topography in a realistic manner.

Whilst the process of forming stencils is quite straightforward in free space, the boundary presents an issue: derivative operators may be truncated, leaving stencil points on the exterior. Forming FD operators in such situations requires the incorporation of additional constraints, in this case derived from the boundary conditions. This concept has precedent in the design of FD operators, aiming to address the lack of information via some other design criterion (e.g. dispersion optimization as in \citet{Tam1993}).

\plot{matrix_intersected_support.png}{width=\columnwidth}
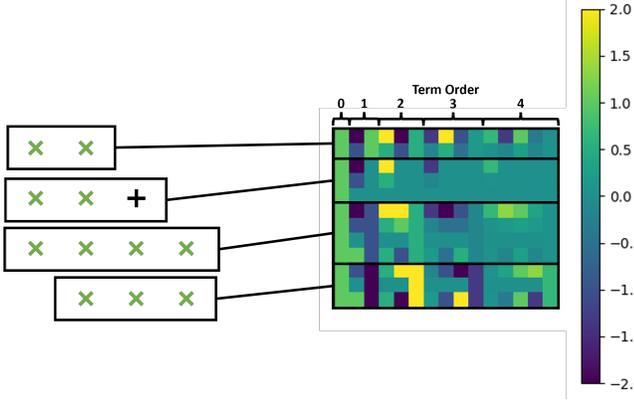
{The effect of truncation of the stencil footprint shown in Figure \ref{fig:matrix_free_space_support.png} by a 45\textdegree  boundary. Note that the points removed from the stencil correspond to rows removed from $\mathbf{A}$.
}

Figure \ref{fig:matrix_intersected_support.png} shows the effect of a boundary cutting through the support region shown in Figure \ref{fig:matrix_free_space_support.png} on the linear system: $\mathbf{A}$ is underdetermined in its current form.

Immersed boundary treatments extend the function to the edge of the domain in a manner that respects the boundary position and conditions, while maintaining consistency with the error introduced by the spatial discretization. This is achieved by including boundary points, which serve to place additional constraints upon the basis, reducing the number of interior points required and imposing the appropriate behaviour. In doing so, additional rows are added to $\mathbf{A}$. The process of forming these rows given boundary conditions and the points at which they are imposed follows in a similar vein to that of forming polynomial approximations of interior function values.

\subsection{Approximating Boundary Conditions}
In 3D, meaningful linear boundary conditions imposed upon a single field (boundary conditions imposed on multiple fields will be discussed in due course) have the general form
\begin{equation}
    \sum_{m=0}^{M_x}\sum_{n=0}^{M_y}\sum_{l=0}^{M_z}\alpha_{mnl}(\mathbf{\mathrm{x}_b}) \frac{\partial^{m+n+l}f}{\partial x^m \partial y^n \partial z^l}(\mathbf{\mathrm{x}_b}) = g(\mathbf{\mathrm{x}_b})
    \label{eq:accmaths:single_boundary_condition}
\end{equation}
where $M_x$, $M_y$, and $M_z$ are as in Equation \ref{eq:accmaths:truncation}. Conditions of this form will not yield the trivial expression $0=g(\mathbf{\mathrm{x}_b})$ when approximated with the basis. The coefficient $\alpha_{mnl}(\mathbf{\mathrm{x}_b})$ may vary with position, although for the applications considered in this paper, these are constant.

Substituting $f$ for its Taylor-series approximation yields
\begin{equation}
\sum_{m=0}^{M_x}\sum_{n=0}^{M_y}\sum_{l=0}^{M_z}\sum_{i=0}^{M_x}\sum_{j=0}^{M_y}\sum_{k=0}^{M_z}\beta_{mnlijk}\frac{\partial^{i+j+k}f}{\partial x^i \partial y^j \partial z^k}(\mathbf{\mathrm{x}}_0) = g(\mathbf{\mathrm{x}_b})\mathrm{,}
\end{equation}
where
\begin{equation}
    \beta_{mnlijk} = \alpha_{mnl}(\mathbf{\mathrm{x}_b})\frac{\partial^{m+n+l}}{\partial x^m \partial y^n \partial z^l}\frac{(x-x_0)^{i}(y-y_0)^j(z-z_0)^k}{i!j!k!} \mathrm{.}
\end{equation}
As is the case with Taylor-series approximations at interior points, this is expressible as a dot product:
\begin{equation}
    \left(\sum_{m=0}^{M_x}\sum_{n=0}^{M_y}\sum_{l=0}^{M_z}\alpha_{mnl}(\mathbf{\mathrm{x}_b})\frac{\partial^{m+n+l}}{\partial x^m \partial y^n \partial z^l}\mathbf{a} \right) \cdot \boldsymbol{\delta} = g(\mathbf{\mathrm{x}_b}) \mathrm{.}
    \label{eq:accmaths:dot_product_bc}
\end{equation}

Note that the derivative vector is identical to that used for Taylor-series expansion at interior points, and thus some set of interior and boundary constraints can be encapsulated as the multiplication of the derivative vector by some matrix of coefficients. There is no distinction between interior and boundary constraints: it is apparent that interior constraints can be formed from Equation \ref{eq:accmaths:single_boundary_condition}, and both can be represented as rows of the linear system. 

\plot{matrix_boundary_support.png}{width=\columnwidth}
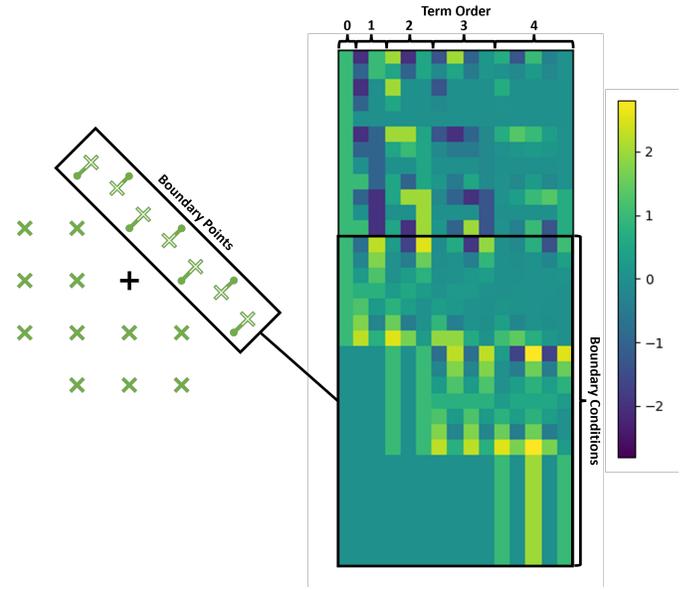
{The effect of adding boundary constraints to the truncated support region shown in Figure \ref{fig:matrix_intersected_support.png}. Boundary points where conditions are imposed are shown as hollow green crosses. Solid green dots show the centre of FD cells containing a boundary point; the normal from this point to the boundary is shown as a green line. Boundary conditions increase in order towards the bottom of the matrix. Note that the highest-order boundary conditions here are invariant with position and thus redundant.
}

With individual constraints for both interior points and boundary conditions, a linear system can be obtained with which to fit the basis, not unlike that in Equation \ref{eq:accmaths:standard_linear_system}. Figure \ref{fig:matrix_boundary_support.png} shows the effect of adding free-surface boundary conditions to the truncated stencil shown in Figure \ref{fig:matrix_intersected_support.png}; the matrix is once again overdetermined, and thus nominally has sufficient information to constrain the derivatives. Note however, that it may be the case that such a matrix is still not full-rank. For example, particular boundary constraints may contain redundant information, as seen in the lowermost rows of the matrix or lack information regarding particular derivatives.

With interior function values and some $J$ boundary conditions imposed at the boundary points $\mathbf{X_b}$, the linear system
\begin{equation}
    \mathbf{A}\boldsymbol{\delta}=
    \begin{pmatrix}
    f(\mathbf{X}) \\
    g_1(\mathbf{X_b}) \\
    \vdots \\
    g_J(\mathbf{X_b})
    \end{pmatrix}\mathrm{,}
    \label{eq:accmaths:general_boundary_linear_system}
\end{equation}
can be formed, with rows of the form given in Equations \ref{eq:accmaths:dot_product} and \ref{eq:accmaths:dot_product_bc}, as shown in Figure \ref{fig:matrix_boundary_support.png}. $\mathbf{X}$ is the set of interior points used to construct the extrapolant and vectors $g_1(\mathbf{X_b})$ through $g_J(\mathbf{X_b})$ contain forcing values corresponding to each boundary condition and forcing point.

As a tangible example, consider constructing a second-order extrapolant in 1D from two interior points (labelled $x_1$ and $x_2$) and a boundary point upon which the boundary conditions $f(x_b)=0$ and $\frac{\partial^2 f}{\partial x^2}(x_b)=0$ are imposed. This particular case yields the linear system
\begin{equation}
    \begin{pmatrix}
        1 & (x_1-x_0) & \frac{(x_1-x_0)^2}{2} \\
        1 & (x_2-x_0) & \frac{(x_2-x_0)^2}{2} \\
        1 & (x_b-x_0) & \frac{(x_b-x_0)^2}{2} \\
        0 & 0 & 1
    \end{pmatrix}
    \begin{pmatrix}
        f(x_0) \\
        \frac{\partial f}{\partial x}(x_0) \\
        \frac{\partial f^2}{\partial x^2}(x_0)
    \end{pmatrix}
    =
    \begin{pmatrix}
        f(x_1) \\
        f(x_2) \\
        0 \\
        0
    \end{pmatrix}
    \label{eq:accmaths:1d_example_system}
    \mathrm{.}
\end{equation}
Whilst higher-order accuracy and higher dimensionality increase the size of this system (making it somewhat unwieldy to show here) the process of constructing this system remains the same. Note here that whilst this particular system is overdetermined, removing a single interior point from this example gives the approach detailed in \citet{Mulder2017} and \citet{Mulder2017a} in which all boundary constraints are selected then supplemented with interior constraints to obtain a square system for which an inverse can be found.

The process of selecting a suitable set of interior and boundary conditions is worthy of some consideration here. Whilst the method used by \citet{Mulder2017, Mulder2017a} yields easy-to-invert systems and is intuitive in 1D, in higher dimensions more choices are available.

\plot{boundary_support.png}{width=0.8\columnwidth}
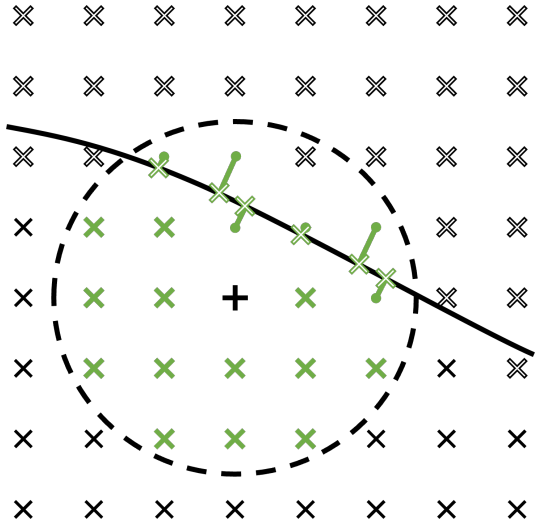
{Construction of the support region in the vicinity of the boundary. Symbols are as in Figures \ref{fig:circular_support.png} and \ref{fig:matrix_boundary_support.png}. The solid black line represents the boundary surface. Hollow black crosses are considered to be outside or too close to the boundary for the purposes of constructing the extrapolant (as defined in the proceeding section).
}

To preserve the Cartesian topology of the grid (thereby allowing all points to be directly indexed), boundary points are taken as the intersection of the boundary normal at a grid point and the boundary, assuming this lies within a point-centred hyperrectangle whose sides correspond to a single grid increment in each dimension. Such a support region is shown in Figure \ref{fig:boundary_support.png}. From here, a cutoff parameter $\eta$ is defined. Interior points where a boundary point lies within a gridpoint-centred hyperrectangle with sides of length $2\eta\Delta x_n$ (where $\Delta x_n$ is the grid increment in the $n$\textsuperscript{th} dimension) are excluded for the purposes of constructing the extrapolant to ensure stability. Note that whilst suitable values of $\eta$ have been empirically determined in this work, more rigorous analysis would be beneficial (although this may be challenging, as in \citealp{Mulder2017a}). An initial support region radius of $(M+1)/2$ is selected, expanded incrementally if insufficient information to constrain the extrapolation is contained therein.

If the basis is sufficiently constrained,  $rank(A)$ will equal the number of expansion terms. The required rank is $M+1$ in 1D, $(M+1)(M+2)/2$ in 2D, and $(M+1)(M+2)(M+3)/6$ in 3D respectively, the latter two corresponding to the $(M+1)$\textsuperscript{th} triangular and tetrahedral numbers respectively.

If one does not intend to be so strict about maintaining formal order, the order of the polynomial basis can be reduced instead, as in \citet{Mulder2017, Mulder2017a}. It is anticipated that this may be advisable with some boundary conditions and geometries to avoid stencil footprints becoming excessively large in edge cases.

\subsection{Reconciliation with the Interior Numerical Discretization}
Whilst the method described above can be used to directly obtain FD stencils, sudden switches in derivative approximation may lead to instability, and it is thus desirable to retain the interior discretization throughout the computational domain. The truncation of interior operators applied in the vicinity of the boundary can be addressed by using the boundary-constrained N-dimensional polynomial basis to approximate function values at required exterior points. The process of projecting this basis onto some set of required exterior points $\mathbf{X_e}$ can be represented by the matrix-vector multiplication
\begin{equation}
    \mathbf{B}\boldsymbol{\delta}
    = \tilde{f}(\mathbf{X_e}) \mathrm{,}
    \label{eq:accmaths:projection_system}
\end{equation}
where $\mathbf{B}$ contains the terms associated with each derivative in the Taylor series evaluated at the respective points in $\mathbf{X_e}$. Derivatives are approximated as
\begin{equation}
    \boldsymbol{\delta} = \mathbf{A}^+
    \begin{pmatrix}
    f(\mathbf{X}) \\
    g_1(\mathbf{X_b}) \\
    \vdots \\
    g_J(\mathbf{X_b})
    \end{pmatrix} \mathrm{,}
\end{equation}
where $\mathbf{A}^+$ is the Moore-Penrose pseudoinverse of $\mathbf{A}$. Note that this is the inverse of Equation \ref{eq:accmaths:general_boundary_linear_system}.

Using the 1D example shown in Equation \ref{eq:accmaths:1d_example_system}, assume two exterior points, designated $x_3$ and $x_4$ are required to complete the stencil operator applied at the specified position. $\mathbf{B}$ will then be of the form
\begin{equation}
    \mathbf{B} =
    \begin{pmatrix}
        1 & (x_3-x_0) & \frac{(x_3-x_0)^2}{2} \\
        1 & (x_4-x_0) & \frac{(x_4-x_0)^2}{2}
    \end{pmatrix} \mathrm{,}
\end{equation}
and thus Equation \ref{eq:accmaths:projection_system} will yield approximations of $f(x_3)$ and $f(x_4)$ as functions of $f(x_1)$, $f(x_2)$ (the other entries in the right-hand side vector of Equation \ref{eq:accmaths:1d_example_system} are zero). More generally, this process yields exterior function values as functions of interior values and boundary conditions.

Substituting these into the interior stencil as necessary, a modified version of the interior operator is obtained. This operator has a support region consisting of the unity of interior points in the original stencil ($\mathbf{X_i}$) and those in the circular stencils used for extrapolation, plus boundary points where conditions are enforced.

A stencil expression approximating some derivative of $f$ can be expressed as
\begin{equation}
    \mathbf{w}\cdot f(\mathbf{X}) = d \mathrm{,}
\end{equation}
where $d$ is some arbitrary derivative and $\mathbf{w}$ contains the stencil weights corresponding to each function value in the vector $f(\mathbf{X})$. Separating the corresponding vector of stencil weights $\mathbf{w}$ into two subvectors $\mathbf{w_i}$ containing weights for points on the interior and $\mathbf{w_e}$ for points on the exterior, the expression for the modified stencil $d$ can be obtained via
\begin{equation}
    \begin{pmatrix}
        \mathbf{w_i} \\
        \mathbf{w_e}
    \end{pmatrix}
    \cdot
    \begin{pmatrix}
        f(\mathbf{X_i}) \\
        \tilde{f}(\mathbf{X_e})
    \end{pmatrix}
    = d\mathrm{.}
\end{equation}
Note that by constructing modified operators in this manner, the extrapolated values are used without being explicitly calculated, with the extrapolation process baked into the stencil coefficients. Furthermore, derivative stencils centred at different gridpoints will have independent extrapolations. This independence and the local nature of these extrapolation operators ensures that the linear systems constructed remain relatively small and simplifies the process of obtaining a solution to the system (both criterion prioritised by \citealp{Hu2016}).

\section*{Application to the 2nd-Order Acoustic Wave Equation}
Application to the 2nd-order acoustic wave equation represents a suitable first step for the detailed method. Containing only a single time-variant field and equation, it has reduced computational cost and implementation complexity compared to other formulations and wave equations; simpler equations also yield simpler boundary conditions. Whilst only P-wave components can be propagated and more complex physics such as anisotropy and viscoelasticity commonly desired for seismic imaging are omitted, applications remain in medical imaging \citep{Guasch2020} and infrasound studies \citep{Kim2014}. Furthermore, it offers a simple platform to test the applicability of the method to various boundary conditions. Doing so paves the way for further diverse boundary conditions introduced with more complex physics.

The equation itself is given as
\begin{equation}
    \frac{\partial^2 p}{\partial t^2} = c^2\nabla^2p + f \mathrm{,}
    \label{2nd_order_acoustic}
\end{equation}
containing a single time-dependent variable $p$, and parameterised via a wavespeed $c$ which varies in space alone. There is also an additional forcing term $f$ which in the context of seismic simulation takes the form of some point source or set thereof.

In seismic applications, this equation is typically discretized with an explicit timestepping scheme, replacing the time derivative with a second-order centred-difference approximation:
\begin{equation}
    \frac{p^{t+1} - 2p^t + p^{t-1}}{\Delta t^2} = c^2\nabla^2p + f \mathrm{.}
\end{equation}
Rearranging for pressure at the forward timestep,
\begin{equation}
    p^{t+1} = 2p^t - p^{t-1} + \Delta t^2c^2\nabla^2p + \Tilde{f}\mathrm{,}
\end{equation}
is obtained, where $\Tilde{f}=f\Delta t^2$. This is the update equation used to estimate the field at the next iteration.

Spatial derivatives are all contained within the Laplacian; for the purposes of an immersed boundary, modified 2nd-derivative operators will need to be generated for each dimension, assuming the boundary is not axially aligned.

Dependent upon the side from which the wave approaches the surface, the reflection coefficient is near 1 or -1, the former when approaching from above whilst the latter whilst approaching from below. We will consider the latter case for now. In this scenario, the topography for all intents and purposes represents an irregular free surface, upon which the condition
\begin{equation}
    p(t, \mathbf{\mathrm{x}_b}) = 0 \mathrm{,}
    \label{pressure_free_suface}
\end{equation}
is to be applied. Using Equation \ref{2nd_order_acoustic}, incrementally higher-order boundary conditions can be derived, these being $\nabla^2p(t, \mathbf{\mathrm{x}_b}) = 0$, $\nabla^4p(t, \mathbf{\mathrm{x}_b}) = 0$, and so forth.

\subsection{Free-Surface Boundary Constraints}
By substituting the polynomial basis into these boundary conditions in the place of pressure, suitable approximations can be formed. For boundary conditions of higher order than the spatial discretization, this will yield the trivial expression $0=0$, and thus such conditions can be discarded.

Considering a fourth-order discretization in 2D, substituting the corresponding basis into the zeroth-order boundary condition yields the equation
\begin{equation}
\begin{aligned}
    &p(x_0, y_0) + (x_b-x_0)\frac{\partial p}{\partial x}(x_0, y_0) + (y_b-y_0)\frac{\partial p}{\partial y}(x_0, y_0) \\
    &+ \frac{(x_b-x_0)^2}{2}\frac{\partial^2 p}{\partial x^2}(x_0, y_0) \\
    &+ (x_b-x_0)(y_b-y_0)\frac{\partial^2 p}{\partial x \partial y}(x_0, y_0) \\
    &+ \frac{(y_b-y_0)^2}{2}\frac{\partial^2 p}{\partial y^2}(x_0, y_0) + \frac{(x_b-x_0)^3}{6}\frac{\partial^3 p}{\partial x^3}(x_0, y_0) \\
    &+ \frac{(x_b-x_0)^2(y_b-y_0)}{2}\frac{\partial^3 p}{\partial x^2 \partial y}(x_0, y_0) \\
    &+ \frac{(x_b-x_0)(y_b-y_0)^2}{2}\frac{\partial^3 p}{\partial x \partial y^2}(x_0, y_0) \\
    &+ \frac{(y_b-y_0)^3}{6}\frac{\partial^3 p}{\partial y^3}(x_0, y_0) + \frac{(x_b-x_0)^4}{24}\frac{\partial^4 p}{\partial x^4}(x_0, y_0) \\
    &+ \frac{(x_b-x_0)^3(y_b-y_0)}{6}\frac{\partial^4 p}{\partial x^3 \partial y}(x_0, y_0) \\
    &+ \frac{(x_b-x_0)^2(y_b-y_0)^2}{4}\frac{\partial^4 p}{\partial x^2 \partial y^2}(x_0, y_0) \\
    &+ \frac{(x_b-x_0)(y_b-y_0)^3}{6}\frac{\partial^4 p}{\partial x \partial y^3}(x_0, y_0) \\
    &+ \frac{(y_b-y_0)^4}{24}\frac{\partial^4 p}{\partial y^4}(x_0, y_0) = 0\mathrm{,}
\end{aligned}
\end{equation}
and similarly, the zero Laplacian and Biharmonic boundary conditions yield
\begin{equation}
\begin{aligned}
    &\frac{\partial^2 p}{\partial x^2}(x_0, y_0) + \frac{\partial^2 p}{\partial y^2}(x_0, y_0) + (x_b-x_0)\frac{\partial^3 p}{\partial x^3}(x_0, y_0) \\
    &+ (y_b-y_0)\frac{\partial^3 p}{\partial x^2 \partial y}(x_0, y_0) + (x_b-x_0)\frac{\partial^3 p}{\partial x \partial y^2}(x_0, y_0) \\
    &+ (y_b-y_0)\frac{\partial^3 p}{\partial y^3}(x_0, y_0) + \frac{(x_b-x_0)^2}{2}\frac{\partial^4 p}{\partial x^4}(x_0, y_0) \\
    &+ (x_b-x_0)(y_b-y_0)\frac{\partial^4 p}{\partial x^3 \partial y}(x_0, y_0) \\
    &+ ((x_b-x_0)^2+(y_b-y_0)^2)\frac{\partial^4 p}{\partial x^2 \partial y^2}(x_0, y_0) \\
    &+ (x_b-x_0)(y_b-y_0)\frac{\partial^4 p}{\partial x \partial y^3}(x_0, y_0) \\
    &+ \frac{(y_b-y_0)^2}{2}\frac{\partial^4 p}{\partial y^4}(x_0, y_0) = 0\mathrm{,}
\end{aligned}
\end{equation}
and
\begin{equation}
    \frac{\partial^4 p}{\partial x^4}(x_0, y_0) + 2\frac{\partial^4 p}{\partial x^2 \partial y^2}(x_0, y_0) + \frac{\partial^4 p}{\partial y^4}(x_0, y_0) = 0\mathrm{,}
\end{equation}
respectively. Each of these equations can be expressed as in Equation \ref{eq:accmaths:dot_product}, with $\boldsymbol{\delta}$ given by
\begin{equation}
\begin{aligned}
    \boldsymbol{\delta}^T = \bigg(&p(x_0, y_0), \frac{\partial p}{\partial x}(x_0, y_0), \frac{\partial p}{\partial y}(x_0, y_0), \\
    &\frac{\partial^2 p}{\partial x \partial y}(x_0, y_0), \frac{\partial^2 p}{\partial y^2}(x_0, y_0), \frac{\partial^3 p}{\partial x^3}(x_0, y_0), \\
    &\frac{\partial^3 p}{\partial x^2 \partial y}(x_0, y_0), \frac{\partial^3 p}{\partial x \partial y^2}(x_0, y_0), \frac{\partial^3 p}{\partial y^3}(x_0, y_0), \\
    &\frac{\partial^4 p}{\partial x^4}(x_0, y_0), \frac{\partial^4 p}{\partial x^3 \partial y}(x_0, y_0), \frac{\partial^4 p}{\partial x^2 \partial y^2}(x_0, y_0), \\
    &\frac{\partial^4 p}{\partial x \partial y^3}(x_0, y_0), \frac{\partial^4 p}{\partial y^4}(x_0, y_0)\bigg) \mathrm{.}
\end{aligned}
\end{equation}
For the zeroth-order condition,
\begin{equation}
\begin{aligned}
    \mathbf{a}^T = \bigg(&1, (x_b-x_0), (y_b-y_0), \frac{(x_b-x_0)^2}{2}, \\
    &(x_b-x_0)(y_b-y_0), \frac{(y_b-y_0)^2}{2}, \frac{(x_b-x_0)^3}{6}, \\
    &\frac{(x_b-x_0)^2(y_b-y_0)}{2}, \frac{(x_b-x_0)(y_b-y_0)^2}{2}, \\
    &\frac{(y_b-y_0)^3}{6}, \frac{(x_b-x_0)^4}{24}, \frac{(x_b-x_0)^3(y_b-y_0)}{6}, \\ 
    &\frac{(x_b-x_0)^2(y_b-y_0)^2}{4}, \frac{(x_b-x_0)(y_b-y_0)^3}{6}, \\
    &\frac{(y_b-y_0)^4}{24}\mathrm{,}
\end{aligned}
\label{eq:bc_vector_0}
\end{equation}
whilst the Laplacian and Biharmonic conditions similarly correspond to
\begin{equation}
\begin{aligned}
    \mathbf{a}^T = \bigg(& 0, 0, 0, 1, 0, 1, (x_b-x_0), (y_b-y_0), (x_b-x_0), \\
    &(y_b-y_0), \frac{(x_b-x_0)^2}{2}, (x_b-x_0)(y_b-y_0), \\
    &(x_b-x_0)^2+(y_b-y_0)^2), (x_b-x_0)(y_b-y_0), \\
    &\frac{(y_b-y_0)^2}{2}\bigg)
\end{aligned}
\label{eq:bc_vector_2}
\end{equation}
and
\begin{equation}
    \mathbf{a}^T = (0, 0, 0, 0, 0, 0, 0, 0, 0, 0, 1, 0, 2, 0, 1)
\label{eq:bc_vector_4}
\end{equation}
respectively. The right-hand side will be zero.

With Equations \ref{eq:bc_vector_0}-\ref{eq:bc_vector_4}, three rows can be constructed for each boundary point used to constrain the basis. Note that the row corresponding with the fourth-order condition is invariant in space, and so will contain the same information irrespective of boundary point location.

\subsection{Application to Example Geometry}
\plot{arc_dy_stencil.png}{width=0.8\columnwidth}
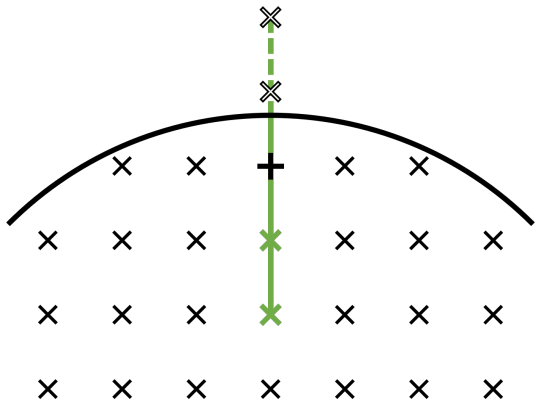
{The stencil footprint of a 4th-order-accurate $\frac{\partial^2 p}{\partial y^2}$ stencil truncated by an arc-shaped boundary. The bold black cross is the stencil position, whilst pale green crosses are other interior points within the stencil. Values at both hollow black crosses are required by the stencil but are located outside the computational domain.
}
To exemplify this process, consider the case of a fourth-order stencil approximating $\frac{\partial^2 p}{\partial y^2}$, truncated by an arc-shaped boundary, as shown in Figure \ref{fig:arc_dy_stencil.png}. It is clear that it will not be possible to form this stencil as one would in free space since two of the requisite points lie outside the physical domain.

\plot{arc_support.png}{width=0.8\columnwidth}
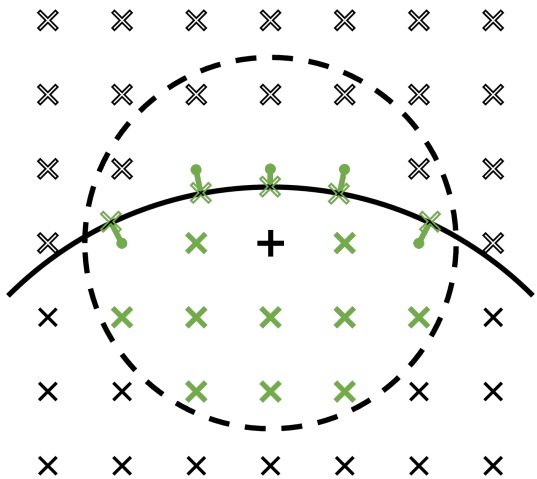
{Support region for a 2D polynomial fitted with a combination of available interior points and boundary points.
}

To rectify this, a circular support region of radius $2.5\Delta x$ is extended from the stencil centre, as in Figure \ref{fig:arc_support.png}. This support radius encircles $5$ boundary points, determined via the previously-discussed criteria. Given a cutoff of $\eta=0.5$ to prevent instability related to the boundary forcing, $11$ interior points are also available for purposes of constructing the 2D extrapolant. Each boundary point will have three matrix rows associated with it, one for each boundary condition imposed at that point, whilst each interior point will correspond to a single row of $\mathbf{A}$.

\multiplot*{2}{arc_support_interior.png,2nd_order_matrix_interior.png}{width=0.6\columnwidth}
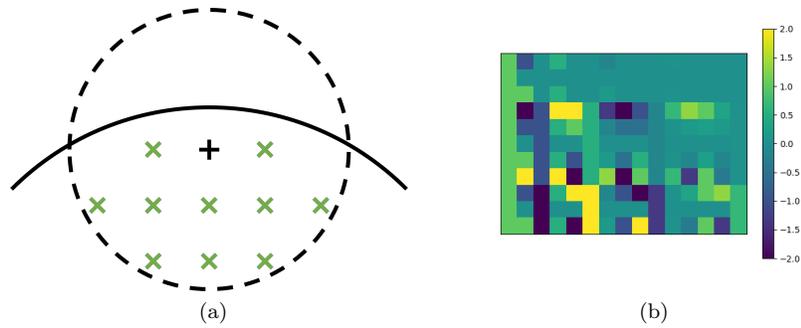
{Interior points within the support region and structure of the corresponding submatrix. Rows in the matrix from top to bottom correspond to points within the support region working downwards from left to right.
}

\multiplot*{2}{arc_support_boundary.png,2nd_order_matrix_multi_points.png}{width=0.6\columnwidth}
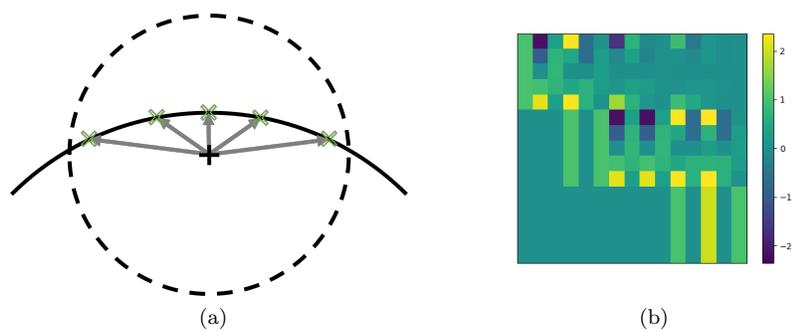
{Boundary points within the support region and structure of the corresponding submatrix. Grey arrows indicate the vector $\mathbf{\mathrm{x}_b} - \mathbf{\mathrm{x}}_0$ for these points. Rows in the matrix from top to bottom correspond to points from left to right.
}

The set of interior points used to fit the basis is shown in Figure \ref{fig:arc_support_interior.png} alongside the resultant submatrix containing the constraints applied at this set of points. Note however that these boundary constraints are not necessarily unique (most prominently the zero biharmonic condition as aforementioned), and may contain overlapping information: whilst the resultant matrix has more columns than rows, it is not necessrily full-rank.

\plot{arc_stencils.png}{width=\columnwidth}
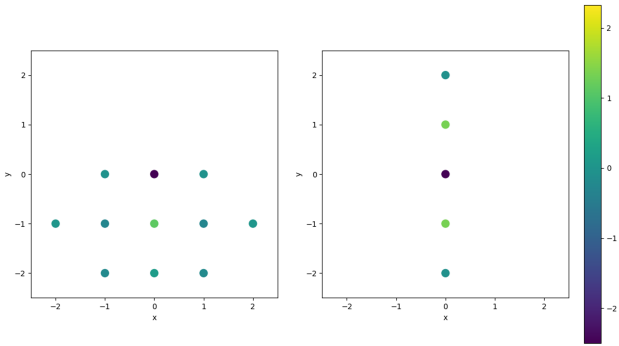
{Comparison of the standard stencil footprint and weights to that of the modified boundary operator for the case illustrated in Figure \ref{fig:arc_dy_stencil.png}. The colourbar indicates the value of the stencil weight at each point.
}

In this particular case, the rank is equal to the number of columns, implying an overdetermined system. To obtain approximations of the derivatives and thus fit the basis, a Moore-Penrose pseudoinverse is used. A weighted least-squares approach was briefly explored to prioritise particular boundary conditions or points but was found to be prone to generating ill-conditioned linear systems whilst having minimal discernible benefit.

Continuing the pressure field onto the pair of required exterior points requires the construction of $\mathbf{B}$ as in Equation \ref{eq:accmaths:projection_system}, evaluating the basis at points $(0, 1)$ and $(0, 2)$, and constructing matrix rows via in the prescribed manner. Note that as all boundary forcing values are zero in this case, corresponding columns of $\mathbf{B}\mathbf{A}^+$ can effectively be ignored for the purpose of constructing the stencil. Applying weights for the interior stencil, the modified boundary operator for this particular point is obtained, as shown in Figure \ref{fig:arc_stencils.png}

\subsection{Reduction to 1D}
Through particular choices made when applying the method described, other, equally feasible stencils can be obtained, including the 1D approximation detailed in \citet{Mulder2017}. Considering the 1D case, the aforementioned free-surface boundary conditions reduce to
\begin{equation}
    p(x_b) = 0, \quad \frac{\partial^2 p}{\partial x^2}(x_b) = 0, \quad \frac{\partial^4 p}{\partial x^4}(x_b) = 0, \quad ...
\end{equation}
and so forth. Suppose some case is encountered where a boundary truncates a stencil as in Figure \ref{fig:arc_dy_stencil.png}. Selecting the $M/2$ closest points to the boundary (at distance greater than $\eta\Delta x$), $\mathbf{A}$ can be constructed such that it is square, enabling derivatives to be obtained by inverting this matrix. Thus the method described by \citet{Mulder2017} can be considered a sub-case of the overarching method described in this paper, albeit distinct from the specific approach taken for the examples shown. The optimal manner by which to delineate the support region and solve the linear system is worthy of future attention - the approach taken by this study is by no means optimal, merely aiming for simplicity and robustness.

\plot{1d_nd_comparison_diff.png}{width=0.95\columnwidth}
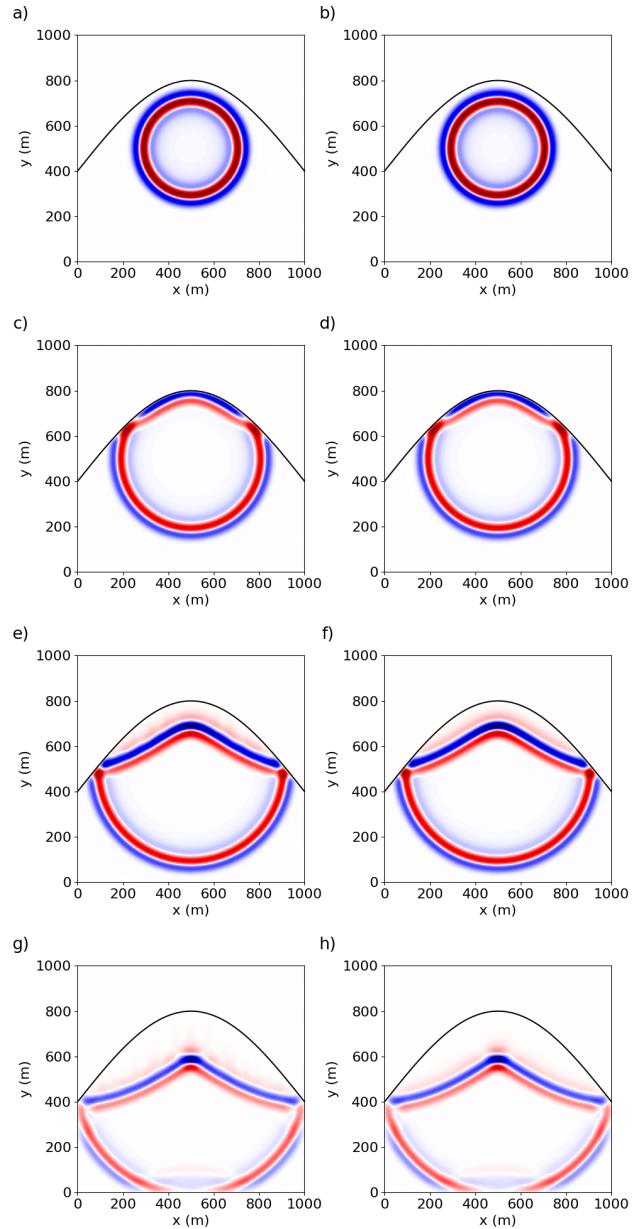
{Snapshots at 400ms, 500ms, 600ms, and 700ms of a wavefront reflecting off a sinusoidal hill upon which a free surface has been imposed. Note that zero Dirichlet boundary conditions have been imposed on all other edges of the domain. Subfigures a, c, e, and g feature an immersed boundary based on a 1D extrapolation scheme, whilst subfigures b, d, f, and h use an N-dimensional basis with circular support. Wavefield amplitudes are normalised against the maximum absolute value in each subfigure for clarity. This convention is continued henceforth.
}

Comparing the proposed approach based on an N-dimensional basis with circular support to an immersed boundary implementation based on per-dimension 1D extrapolations for a 2D second-order acoustic example featuring a sinusoidal free surface, shown in Figure \ref{fig:1d_nd_comparison_diff.png}, both approaches yield visually similar results. However, some minor unevenness is observable in the trailing edge of the reflected wave when 1D extrapolations are used, this area appearing smoother when N-dimensional extrapolation is used.

\subsection{Convergence testing}

\multiplot*{2}{exact_solution.png,convergence_plot.png}{width=0.8\columnwidth}
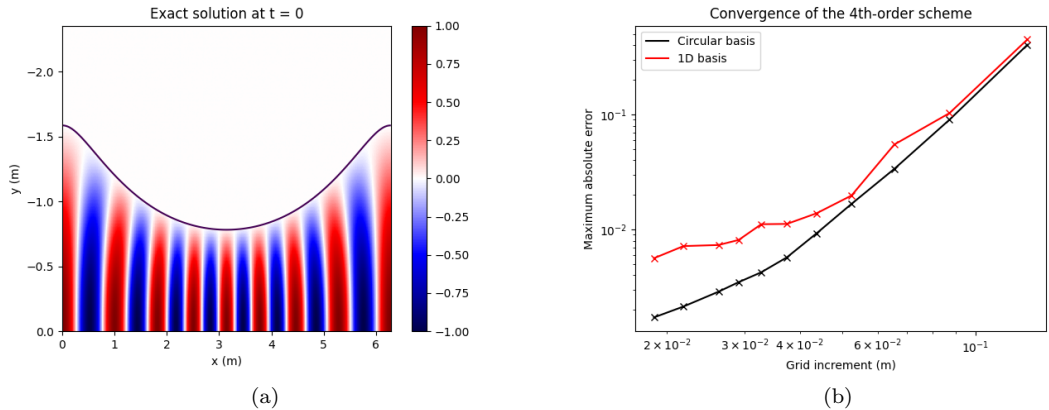
{The exact solution at $t=0$ is shown in the left subfigure. The solid black line is the free surface, the left and right sides of the domain have periodic boundary conditions applied, and the wavefield is mirrored across the lower boundary. A comparison of the convergence behaviour of the method proposed to one using 1D approximations is shown on the right.
}
To examine the convergence behaviour of the proposed approach, a setup initially presented by \citet{Mulder2017} and subsequently used in \citet{Mulder2017a} (the second example in both) was replicated, as shown in Figure \ref{fig:exact_solution.png}. An exact solution exists for this example, allowing for the error in any numerical solution to be evaluated. Figure \ref{fig:convergence_plot.png} shows the convergence behaviour of a scheme based on N-dimensional extrapolations, compared against a scheme based on axially-aligned 1D extrapolations.

For a fourth-order spatial discretization, the reduction in observed error with respect to grid increment was found to be initially just short of fourth-order, flattening around a grid increment of $0.02$ as the spatial error is eclipsed by second-order timestepping error for finer grids. As the timestep was set at 10\% of the critical value, this implies that the error introduced by the immersed boundary was minimal in all cases and that in many cases, topography implementation will cease to be the accuracy bottleneck when this immersed boundary treatment is applied. Reducing the timestep enabled the continuation of the approximate fourth-order trend to finer grids, but accumulation of floating-point error again resulted in a similar flattening albeit at a smaller grid increment. The immersed boundary approach based on N-dimensional extrapolations was found to yield reduced error versus that based on 1D extrapolations for all grid increments tested, particularly at finer resolutions, albeit with similar convergence behaviour.

\section{Extension to Multiple Fields and the 1st-Order Acoustic Wave Equation}
The method detailed yields accurate results for equations concerning a single field, whilst enabling higher-dimensional conditions to be imposed on the boundary. Another appeal of this approach is the readiness with which it is extended to cases where multiple fields are present. The acoustic wave equation can alternatively be formulated as a coupled system of pressure and particle velocity, introducing a spatially-variant density parameter, capturing density-dependent amplitude variations. These equations take the form
\begin{equation}
    \frac{\partial p}{\partial t} = \rho c^2 \boldsymbol{\nabla}\cdot\mathbf{v} + f, \quad \frac{\partial \mathbf{v}}{\partial t} = \frac{1}{\rho}\boldsymbol{\nabla}p\mathrm{,}
    \label{eq:1st_order_wave_equation}
\end{equation}
where $\mathbf{v}$ is particle velocity and $\rho$ is density. This formulation introduces additional fields in the form of components of the particle velocity vector. Given the aforementioned free-surface condition imposed on the pressure field, it is apparent from Equation \ref{eq:1st_order_wave_equation} that the condition
\begin{equation}
    \boldsymbol{\nabla}\cdot \mathbf{v}(t, \mathbf{\mathrm{x}_b}) = 0
\end{equation}
must also be imposed. Whilst boundary conditions considered prior to now have concerned some property of a single field at the boundary, when multiple fields are present in a model, boundary conditions specifying some relationship between these may be present. Each vector component can be approximated by an independent polynomial basis in free space, but at the edge of the domain, they will require construction such that this relationship is respected. As will be highlighted - this extension can be naturally handled by the method described.

The proceeding section considers boundary conditions of this type more generally, before honing in on the application of this approach to the particle velocity free surface.

\subsection{Boundary Conditions Spanning Multiple \\Fields}
 Where multiple fields are present, individual Taylor Series are used to approximate each. Supposing $K$ separate fields, labelled $f_1$ through $f_K$ are present within a model: ignoring boundary conditions for now, the polynomial fitting process can be expressed as
\begin{equation}
    \begin{pmatrix}
        \mathbf{A_1} & & \\
        & \ddots & \\
        & & \mathbf{A_K}
    \end{pmatrix}
    \begin{pmatrix}
        \boldsymbol{\delta_1} \\
        \vdots \\
        \boldsymbol{\delta_K}
    \end{pmatrix}
    =
    \begin{pmatrix}
        f_1(\mathbf{X_1}) \\
        \vdots \\
        f_K(\mathbf{X_K})
    \end{pmatrix}
    \label{eq:accmaths:multiple_fields_system}
    \mathrm{,}
\end{equation}
where $\mathbf{A_k}$ is the matrix containing coefficients associated with each derivative in the polynomial expansion approximating field $f_k$ at the set of points $\mathbf{X_k}$. The vector $\boldsymbol{\delta_k}$ contains the various derivatives of $f_k$, analogous to the single-field case. The left-hand matrix is block-diagonal, and thus the linear system can be split into $K$ smaller systems, each to be solved individually. In the case that boundary conditions concern only a single field apiece, this remains true, and one can still separate the system in this manner. The intuitive implication of this is the independence of polynomial expansions approximating each field in the absence of any constraints which would otherwise link them.

However, if boundary conditions impose some particular relationship between fields, then it is apparent that the resultant polynomial approximations of one of these fields will require information regarding the other fields present in the boundary condition for it to be respected. As in the single-field case, each boundary condition will be approximated with Taylor-series expansions, although note that in this case, each field present within the constraint will have its own expansion. Each of these series contains derivatives of its respective function, contained within the corresponding $\boldsymbol{\delta_k}$ as in Equation \ref{eq:accmaths:multiple_fields_system}. As such, the Taylor-series approximation of a constraint specifying some relationship between multiple fields will contain derivatives of multiple fields, and thus when represented in the previously-detailed dot-product form, the vector $\boldsymbol{\delta}$ will consist of multiple $\boldsymbol{\delta_k}$.

In the context of the linear system shown in Equation \ref{eq:accmaths:multiple_fields_system}, rows corresponding to such boundary conditions will span multiple previously-separate blocks, thereby linking them accordingly. As the process of fitting the extrapolant relies on the inversion of the left-hand-side matrix, it follows that these linked blocks will need to be inverted in tandem as they are no longer separated. As a boundary condition row in such a scenario maps multiple $\boldsymbol{\delta_k}$ onto a single boundary forcing value, and each of the previously-separate blocks maps its respective $\boldsymbol{\delta_k}$ onto the corresponding $f_k(\mathbf{X_k})$, the inverse of the resultant block will consequently map all $f_k(\mathbf{X_k})$ present onto any given derivative.

As a simple example, consider some pair of 1D functions $f$ and $h$ upon which the condition
\begin{equation}
    f(x_b)+h(x_b) = 0\mathrm{,}
\end{equation}
is to be applied at the boundary. Approximating this constraint with 2nd-order Taylor series expanded around some $x_0$ yields
\begin{equation}
\begin{aligned}
    &f(x_b) + (x_b-x_0)\frac{\partial f}{\partial x}(x_b) + \frac{(x_b-x0)^2}{2}\frac{\partial^2 f}{\partial x^2}(x_b) \\
    &+ h(x_b) + (x_b-x_0)\frac{\partial h}{\partial x}(x_b) \\
    &+ \frac{(x_b-x0)^2}{2}\frac{\partial^2 h}{\partial x^2}(x_b) = 0\mathrm{.}
\end{aligned}
\end{equation}
Representing this as a dot product of two vectors
\begin{equation}
    \begin{pmatrix}
        \mathbf{a_f} \\
        \mathbf{a_h}
    \end{pmatrix}
    \cdot
    \begin{pmatrix}
        \boldsymbol{\delta_f} \\
        \boldsymbol{\delta_h}
    \end{pmatrix}
    =0
\end{equation}
where
\begin{equation}
    \mathbf{a_f}^T = \mathbf{a_h}^T = \left(1, (x_b-x_0), \frac{(x_b-x0)^2}{2} \right)\mathrm{,}
\end{equation}
\begin{equation}
    \boldsymbol{\delta_f}^T = \left(f(x_b), \frac{\partial f}{\partial x}(x_b), \frac{\partial^2 f}{\partial x^2}\right)\mathrm{,}
\end{equation}
and
\begin{equation}
    \boldsymbol{\delta_h}^T = \left(h(x_b), \frac{\partial h}{\partial x}(x_b), \frac{\partial^2 h}{\partial x^2}\right)\mathrm{.}
\end{equation}
To avoid confusion, note that $\mathbf{a_f} = \mathbf{a_h}$ results from the boundary condition applied and will generally not be the case, depending upon on constraints applied. A single point of $f$ in this form is given as
\begin{equation}
    \begin{pmatrix}
        \mathbf{a_f} \\
        0
    \end{pmatrix}
    \cdot
    \begin{pmatrix}
        \boldsymbol{\delta_f} \\
        \boldsymbol{\delta_h}
    \end{pmatrix}
    =0\mathrm{.}
\end{equation}
Approximations of $h$ at interior points can be formed in a similar manner. Assembling values of $f$ and $h$ at some set of points into a linear system using the above methodology, it is apparent that the row corresponding to the boundary condition links two otherwise separate blocks pertaining to interior points of $f$ and $h$. A more thorough discussion on the effect of such boundary conditions on the structure of the linear system can be found in Appendix \ref{matrixstructure}.

Solving the resultant linear system gives derivatives of the basis (now multiple bases) at the expansion point as in the single-field case. Given these, the function can be continued beyond the boundary to obtain values at exterior points required by interior stencils. This is achieved in the same manner as the single-field case. These function values will be some weighted sum of boundary forcing values and interior values of all fields linked to the function of interest via boundary constraints. Consequently, any stencil formed using these approximated values will span all of these fields as well.

With this approach, the need to devise application-specific approximations is removed, enabling the application of a consistent method across a wide range of boundary conditions. As such, for some given set of derivative operators and boundary conditions (within the parameters discussed above), a scheme of this class can be generated.

\subsection*{The Particle Velocity Free-Surface}
Returning to the zero-divergence boundary condition to be imposed on the particle velocity vector, approximating this with a second-order basis gives
\begin{equation}
\begin{aligned}
    &\frac{\partial v_x}{\partial x}(x_0, y_0) + (x_b-x_0)\frac{\partial^2 v_x}{\partial x^2}(x_0, y_0) \\
    &+ (y_b-y_0)\frac{\partial^2 v_x}{\partial x \partial y} + \frac{\partial v_y}{\partial y}(x_0, y_0) \\
    &+ (x_b-x_0)\frac{\partial^2 v_y}{\partial x \partial y} + (y_b-y_0)\frac{\partial^2 v_y}{\partial y^2}(x_0, y_0) = 0 \mathrm{,}
\end{aligned}
\end{equation}
yielding the characteristic row of $\mathbf{A}$:
\begin{equation}
\begin{aligned}
    \mathbf{a}^T = \bigg(&0, 1, 0, (x_b-x_0), (y_b-y_0), 0, \\
    & 0, 0, 1, 0, (x_b-x_0), (y_b-y_0)\bigg) \mathrm{.}
\end{aligned}
\end{equation}
There will be one such row for every boundary point. Returning to the geometry shown in Figure \ref{fig:arc_support.png}, boundary conditions are imposed at a single set of points for all fields to ensure consistency, and as such the particle velocity boundary condition will be imposed on these same points. Whilst the support region will nominally be of a smaller radius in this case and only contain the middle three boundary points, it happens that for this particular geometry, there is insufficient information within this radius, and thus it must be expanded to that shown in Figure \ref{fig:arc_support.png} to constrain the basis.

\plot{1st_order_matrix_assembly.png}{width=\columnwidth}
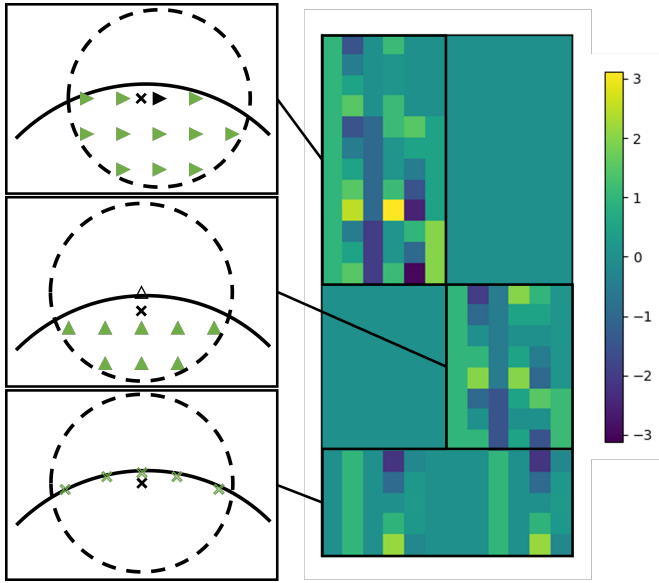
{Assembling $\mathbf{A}$ for the particle velocity fields. Horizontal arrows correspond with $v_x$ points used to construct the extrapolant, whilst vertical arrows are the respective $v_y$ points. The black cross designates the expansion point. Note that only the subgrids of concern are shown for clarity. Note how the interior points of each field correspond with an independent block, whilst the boundary condition rows span both.
}

To prevent errors and nonphysical effects introduced by overextended velocity field extrapolations, $\eta$ is set to zero for both velocity fields. As such, the support region for each particle velocity component is as shown in Figure \ref{fig:1st_order_matrix_assembly.png} alongside the corresponding $\mathbf{A}$. The structure of $\mathbf{A}^+$ is shown in Figure \ref{fig:1st_order_matrix_inverse.png}:
\plot{1st_order_matrix_inverse.png}{width=0.8\columnwidth}
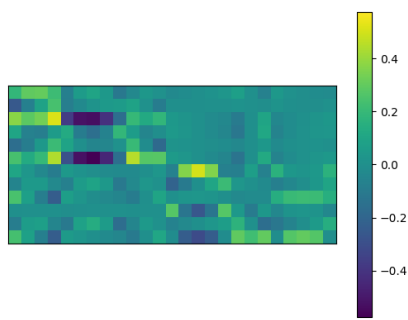
{The structure of the Moore-Penrose pseudoinverse of $\mathbf{A}$ shown in Figure \ref{fig:1st_order_matrix_assembly.png}.
}
it is apparent that this matrix is dense. Given it maps values of both velocity components and boundary forcing values onto derivatives of both basis functions, each of these derivatives will be approximated as a weighted sum of values of both components. Consequently, any extrapolations using these bases will also be in terms of both $v_x$ and $v_y$ and thus stencils completed using these extrapolations will span both these fields too.

Note that a staggered scheme is used to prevent the emergence of checkerboard instability when solving this particular formulation. As such the base stencil used to construct the modified FD operators will be backward-staggered.

The necessity of expanding the support region in this case is worth revisiting. The particle velocity free-surface and accompanying higher-order conditions
$\nabla^2\boldsymbol{\nabla}\cdot \mathbf{v}(\mathbf{\mathrm{x}_b}) = 0$, $\nabla^4\boldsymbol{\nabla}\cdot \mathbf{v}(\mathbf{\mathrm{x}_b}) = 0$, and so forth concern multiple fields. Consequently, the number of boundary constraints within a support region of some given radius is likely to be low relative to the number of fields, requiring an enlarged support region in some cases.

\subsection{Convergence Testing}
To explore the convergence behaviour of the boundary treatment devised, the previous setup was replicated for testing with the 1st-order acoustic wave equation. As before, a 4th-order spatial discretization was used and the timestep was set to 10\% of the critical timestep.

\plot{convergence_plot_1st_order.png}{width=\columnwidth}
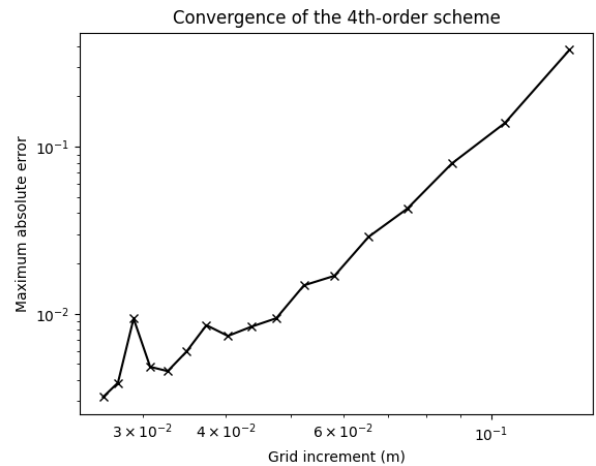
{Convergence of the numerical scheme for the 1st-order formulation of the acoustic wave equation using the same setup shown in Figure \ref{fig:exact_solution.png}. Note that scales are slightly different in this figure.
}

As in the previous test, initial convergence is just short of 4th-order, before gradually flattening around a grid increment of $0.01$, at which point convergence is around 2nd-order as timestepping error saturates the solution. Note that convergence is somewhat less smooth beyond this point with a handful of blackspots where error is anomalously high compared to the prevailing trend. However, broadly speaking, the maximum error in the scheme continues to fall as the discretization is refined. Investigation of these anomalies found them to be specific to very particular grid sizes (adding or removing a single node is sufficient to prevent the more prominent spike), although the reason for this remians to be determined. Overall however, it appears that the error introduced by the boundary is rapidly eclipsed by other sources as the discretization is refined, implying that it is unlikely to be a significant source of error in practical applications.

\section*{Implementation}
 Given a set of symbolic equations that hold on the boundary surface and a discretized signed-distance function (SDF) encapsulating the boundary position, a suitable numerical scheme and thus modified stencils implementing the immersed boundary can be automatically devised. A framework to do so, Schism, was developed as a plugin for Devito. This was done not only to expediate and simplify the implementation of the following test cases but to explore synergies between this generalised immersed boundary method and code generation. Due to the high-level nature of the abstractions created, the introduction of an immersed boundary to a numerical model written in Devito can be achieved with only a handful of additional lines of code and a qualitative understanding of what is being done behind the scenes.

This enabled all the examples shown in this paper to be implemented with a common codebase - only the top-level model specification is changed between examples. The unprecedented flexibility of this approach enables a wide range of geophysical scenarios to be condensed into an understandable, repurposable form, enabling maximum code reuse. Whilst a comprehensive overview of the mechanisms by which this was achieved is beyond the scope of this paper, the proceeding examples all leverage this functionality.

\section*{Examples}
\renewcommand{\figdir}{Examples} 
Reflecting the wide range of relevant geophysical applications, we present a suite of examples showcasing our approach. These examples are designed to resemble particular problems of interest and are based on real-world topography.

\subsection*{2D 2nd-order acoustic free-surface}
The first example is based on an East-West profile across the summit of Mount St Helens, Washington, USA. The summit collapse during the 1980 eruption and subsequent lava dome formation within the crater resulted in near-vertical crater walls and a mixture of concavity and convexity on the crater floor. As a stratovolcano, Mount St Helens has steep, uneven flanks rising over a kilometre from the surrounding landscape, making it an ideal stress test for the method proposed.

The topographic surface was represented internally as an SDF discretised onto the FD grid. This representation has several advantages; its mathematical properties lend themselves to straightforward geometry handling, and it ensures that the resolution of the surface and resultant accuracy of the surface representation within the numerical scheme are consistent with the interior discretization. Note however that this representation enables the surface to be located with much finer precision than the FD discretization, despite its discretization on the same grids, and the SDF can be constructed from extremely high-resolution digital elevation models (DEMs) without the requirement to downsample or alias the raw data to match the FD grid. This is beneficial for real-world applications where such data obtained from satellites and drones may be structured or unstructured (or compounds of multiple such datasets), heavily oversampled versus the discretization required for numerical-dispersion-free wavefield propagation, and typically with a vertical precision in the order of centimetres.

Material properties are kept consistent throughout the model domain, such that all topographic interactions observed are a product of the boundary treatment rather than any material discontinuity: a convention continued throughout the examples presented. Variable material parameters are a separable concern and are essentially trivial to implement (particularly with the abstraction layers used in this study). Furthermore, their introduction runs the risk of inconsistency between implicit interfaces represented by material parameters and the explicit interface encapsulated by the immersed boundary. On this basis, it is not recommended to include any material contrast at the surface.

As aforementioned, a free-surface boundary condition is imposed on the pressure field on the topographic surface. The remaining edges of the computational domain have zero Dirichlet boundary conditions imposed for simplicity. In practice, one would apply a damping boundary condition of choice along these edges, but again this is a separable concern and straightforwardly combined with the method presented in this paper.

A Ricker wavelet with a peak frequency of 8Hz is injected below the lava dome, at an elevation of 1250m, in a location loosely reminiscent of tremors induced by shallow magma movement within the conduit. Placing the source close to the surface maximises observed interaction between wavefield and topography. $\Delta x=\Delta y=30m$ and the Courant number is set to $0.5$. The spatial discretization used is fourth-order accurate: a precedent continued to all other examples shown.

\plot{2d_free_surface_helens.png}{width=\columnwidth}
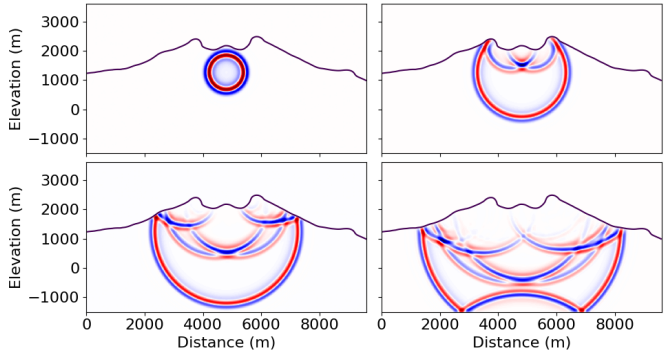
{Snapshots of the wavefield interacting with a free surface at 375ms, 750ms, 1125ms, and 1500ms. Celerity is 2.5km/s throughout the model. The black line designates the isosurface $s(x, y)=0$ on the SDF, coinciding with the surface. The wavefield shown in each snapshot is normalised for clarity.
}

We see in Figure \ref{fig:2d_free_surface_helens.png} that the uneven topography results in several distinct reflections, with further minor reverberations trailing the main wavefront. Also apparent is the focusing and defocusing effect of concave and convex topography respectively, and the diffraction of the wavefront around obstructions. This yields a much more complex series of arrivals than would be observed for a flat surface, although the horizontally-propagating wavefront is only mildly distorted, in agreement with previously-published findings.

\subsection*{2D 2nd-order acoustic rigid-surface}
The free surface is not the only boundary condition of geophysical interest in the context of wavefield propagation. For acoustic waves propagating in the air, the Earth represents an extremely dense and essentially immobile surface: zero particle velocity at the interface corresponds to a rigid surface. The Mount St Helens profile used in the previous example once again features, although the surface now forms the lower bound of the domain.

In an imitation of a typical infrasound propagation scenario, a 1Hz Ricker source is placed at an elevation of 2600m above the lava dome. Located only slightly above the crater rim, this location was chosen to capture reverberation within the crater without trapping the majority of the wave within.

\plot{2d_rigid_surface_helens.png}{width=\columnwidth}
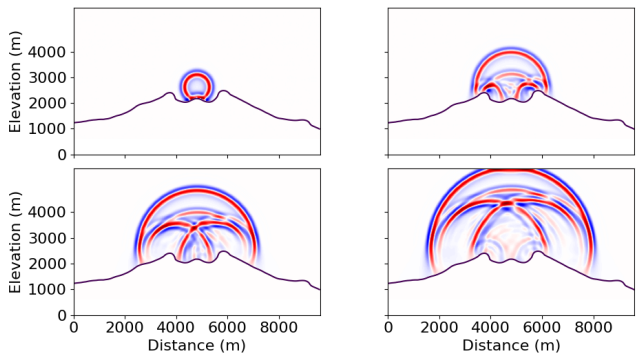
{Snapshots of the wavefield interacting with a rigid surface at 2.5s, 5s, 7.5s, and 10s. Celerity is 350m/s throughout the model. Parameters are altered in this example to better reflect infrasound propagation problems to which this boundary condition is applicable. The wavefield shown in each snapshot is normalised for clarity.
}

Figure \ref{fig:2d_rigid_surface_helens.png} shows the propagating wavefield, including the multiple distinct reflected wavefronts. Note the reversed polarity of these wavefronts versus those found in the previous example. Considerable reverberation within the crater is observed, alongside the wavefront diffracting over obstacles (most notably the crater rim).

\subsection*{3D 2nd-order acoustic free-surface}
Whilst in some cases, wave propagation within a 3D physical domain can be adequately approximated along some suitable transect, this generally relies on some kind of consistency within the omitted dimension. Such approximations are possibly less suitable in the setting of a volcanic edifice, containing strong topographic variation along all directions: particularly true for Mount St Helens due to the collapsed northern flank. In such settings, full 3D modelling may be necessary to achieve realistic wavefield propagation.

\plot{3d_free_surface_render.png}{width=\columnwidth}
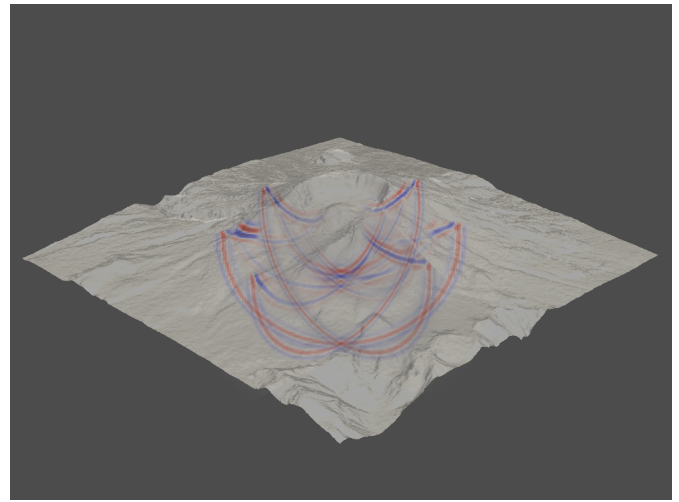
{Render of the 3D free-surface wavefield at 1125ms and topography. Slices of the wavefield are shown aligned and diagonal to each compass direction for clarity. Wavefield transparency is scaled with amplitude to emphasise the wavefronts.
}
The setup (besides the obvious) for this example was much the same as the prior 2D free-surface example. Figure \ref{fig:3d_free_surface_render.png} shows the results of wavefield interaction with the volcanic topography. The first arrivals radiate outwards with minimal interruption, diffracting around obstacles in their path, with complex frills of reflected wavefronts from larger-scale surfaces forming a layered series of  distinct arrivals, leaving minor reverberations in their wake. With only the first arrivals remaining relatively unscathed, this illustrates the error in assuming that a flat surface will adequately reproduce wavefield behaviour observed in rough terrain.

\plot{3d_free_surface.png}{width=\columnwidth}
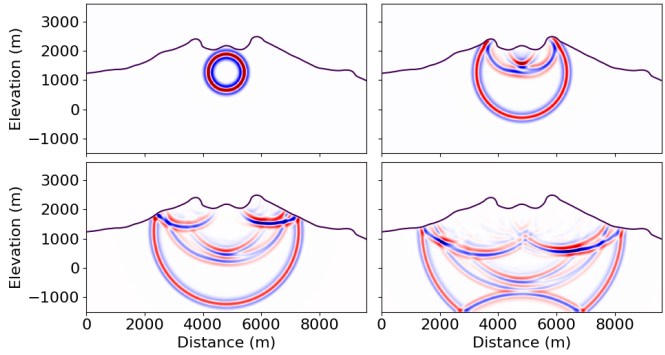
{Slices through the wavefield interacting with the 3D free-surface on the profile used for the 2D examples (the x-z plane). Snapshots were taken at 375ms, 750ms, 1125ms, and 1500ms respectively. Celerity is 2.5km/s throughout the model. The black line designates the isosurface $s(x, y, z)=0$ on the SDF, coinciding with the surface. The wavefield shown in each snapshot is normalised for clarity.
}

Wavefield slices shown in Figure \ref{fig:3d_free_surface.png} are at first glance similar to those of their 2D approximation in Figure \ref{fig:2d_free_surface_helens.png}: the position, shape, and relative amplitude of major wavefronts are the same. On closer inspection however, greater complexity emerges: the main wavefronts contain additional reflections and smaller-scale differences are visible, particularly around the crater rim where some out-of-plane reverberations make themselves known.

\subsection*{3D 2nd-order acoustic rigid-surface}

\plot{3d_rigid_surface_render.png}{width=\columnwidth}
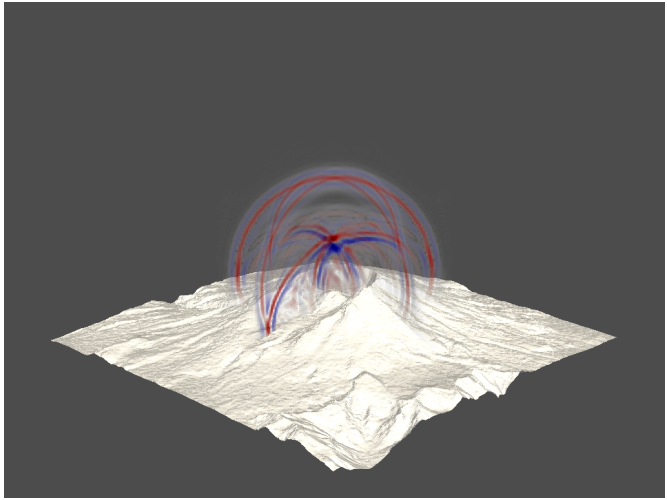
{Render of the 3D rigid-surface wavefield at 7.5s and topography. Slices of the wavefield are shown aligned and diagonal to each compass direction for clarity. Wavefield transparency is scaled with amplitude as before.
}

The previous rigid-surface example is similarly extended to 3D, the results of which are shown in Figure \ref{fig:3d_rigid_surface_render.png}. In this case, likely due to the concavity of the geometry in the vicinity of the source, even more pronounced out-of-plane reflections are observed, particularly within the confines of the caldera where complex and protracted reverberations can be clearly seen. Again, the most prominent reflected wavefronts become more confused in 3D, exhibiting less coherency and continuity due to the wide range of paths taken by its constituent reflections.

\plot{3d_rigid_surface.png}{width=\columnwidth}
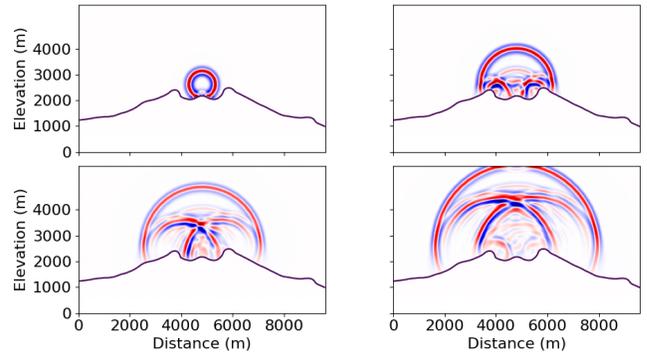
{Snapshots of the wavefield interacting with the 3D rigid surface at 2.5s, 5s, 7.5s, and 10s. The transect is chosen to match that used for the 2D examples (slices on the x-z plane). Celerity is 350m/s throughout the model. The wavefield shown in each snapshot is normalised for clarity.
}

These effects are particularly apparent when comparing Figure \ref{fig:3d_rigid_surface.png} to Figure \ref{fig:2d_rigid_surface_helens.png}, demonstrating the strongly 3D nature of wavefield propagation in this scenario.

\subsection*{2D 1st-order acoustic free-surface}
As aforementioned, the detailed immersed boundary approach is equally applicable to the 1st-order formulation of the acoustic wave equation, and to this end, the setup from the prior 2D 2nd-order acoustic free-surface example is revisited, the results of which are shown in Figure \ref{fig:2d_1st_order_helens.png}.

\plot{2d_1st_order_helens.png}{width=\columnwidth}
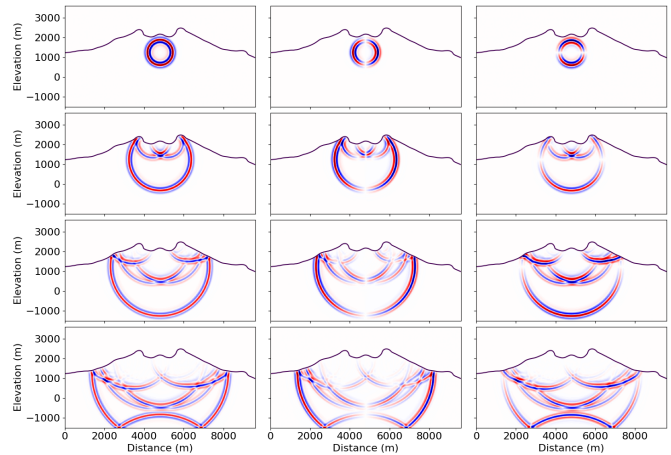
{Snapshots of the pressure(left), x particle velocity (middle), and y particle velocity (right) wavefields at 375ms, 750ms, 1125ms, and 1500ms. Celerity is 2.5km/s throughout the model, density is homogeneous throughout. The black line designates the isosurface $s(x, y)=0$ on the SDF, coinciding with the surface. The wavefields shown in each snapshot are normalised for clarity.
}

The same reflection geometry as in the 2nd-order example is observed, with the additional particle velocity highlighting the partitioning of energy between horizontal and vertical particle motion, dependent on the orientation of the reflector. The success of the vector boundary condition implementation is clear in the clean, artefact-free reflections observed.

\subsection*{3D 1st-order acoustic free-surface}
\plot{3d_free_surface_render_1st_order.png}{width=\columnwidth}
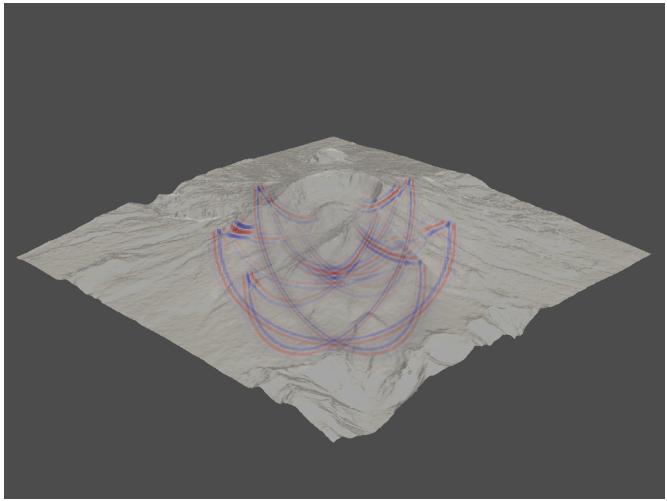
{Render of the 3D pressure wavefield at 1125ms and topography. Slices of the wavefield are shown aligned and diagonal to each compass direction for clarity. Wavefield transparency is scaled with amplitude to emphasise the wavefronts.
}

As is precedent at this point, the previous 2D model is extensible to 3D to demonstrate the capabilities of this approach, once again reusing the prior parameterisations.

\plot{3d_free_surface_1st_order.png}{width=\columnwidth}
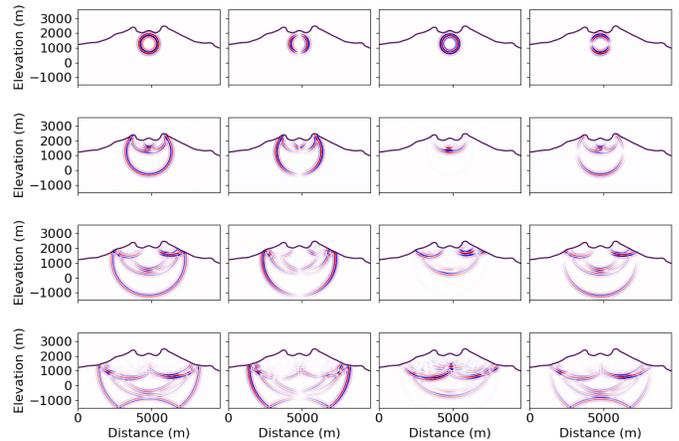
{Snapshots of the pressure (left), x particle velocity (middle left), y particle velocity (middle right), and z particle velocity (right) wavefields at 375ms, 750ms, 1125ms, and 1500ms. The y-axis is oriented into the page.  Celerity is 2.5km/s throughout the model, density is homogeneous throughout. The black line designates the isosurface $s(x, y, z)=0$ on the SDF, coinciding with the surface. The wavefields shown in each snapshot are normalised for clarity.
}

From Figure \ref{fig:3d_free_surface_render_1st_order.png}, it is apparent that the pattern of radiation is identical the that with the 2nd-order formulation (note that the wavelet shape changes between the two, as the source time series was kept the same). Figure \ref{fig:3d_free_surface_1st_order.png} shows snapshots of the wavefields as the wave propagates. The y particle velocity field exemplifies the strongly 3D nature of interaction with topography, with reflected energy into the page strongly apparent. As with the 2nd-order formulation, reflections become more complex and confused in 3D due to the nature of the topography.

For this run, it was found that particle velocity stencils did become somewhat larger than would be desirable at a small handful of points, highlighting limitations of the Taylor series as a basis. It is anticipated that with an improved choice of basis (and potentially linear system setup and solver), the support region could be reduced to a more manageable size. Alternatively, the construction of a more targeted support region may aid in alleviating this issue. However, for this particular run, a basis-order-reduction strategy (as used by \citealp{Mulder2017}) was used at points where insufficient information was present to construct the extrapolant in an effort to rein in this stencil growth.

\section*{Conclusions}
A general immersed boundary treatment is presented, enabling a consistent methodology to be applied across multiple wave equations and boundary conditions. The boundary is encapsulated by modified FD stencils with a circular support region and spatially-variant coefficients, using an N-dimensional Taylor-series extrapolation scheme to continue the field beyond the edge of the domain. As the approach proposed naturally accommodates the implementation of higher-dimensional and vector boundary conditions, it is not necessary to make any application-specific approximations to the boundary. The efficacy of this approach was demonstrated via convergence tests and a range of numerical examples featuring real-world topography, implementing both free and rigid surfaces with the first and second-order acoustic wave equations in 2D and 3D. The one-size-fits-all nature of the method presented enabled the development of a high-level framework, Schism, allowing all of these examples to be implemented via a broadly common codebase. This approach to immersed boundary implementation synergises with emerging code-generation approaches to FD kernel implementation, which was leveraged throughout this paper.

\section{Code Availability and Reproducability}
Schism is an open-source codebase and can be found at github.com/devitocodes/schism. This repository contains all the examples shown in this paper (alongside others), and a suite of unit tests. Schism can be installed from this repository as a Python module using Pip, or alternatively, a Dockerfile to run the code is also included. The codebase at the time of publication is archived on Zenodo at https://zenodo.org/record/8167794.

\section{Acknowledgements}
We wish to thank Wim Mulder for providing assistance with replicating his exact solutions, and Tim MacArthur for exploring the capabilities of our codebase in his own experiments. We also extend thanks to the rest of the Devito team and wider community for their feedback and support, without which this work would not have been possible. This work was funded as part of EPSRC DTP training grant EP/R513052/1.

\append[matrixstructure]{Boundary Effect On Matrix Structure}
\renewcommand{\figdir}{Appendices} 
The boundary serves both to truncate the support region (as no data lies beyond it) and to introduce additional constraints to the polynomial fitting. This profoundly alters the linear system that must be solved to fit the polynomial, particularly where boundary conditions linking fields are present: in the case of the particle velocity boundary condition for the free surface for example.

\multiplot*{2}{matrix_support_no_boundary.png,matrix_no_boundary_annotated.png}{width=0.8\columnwidth}
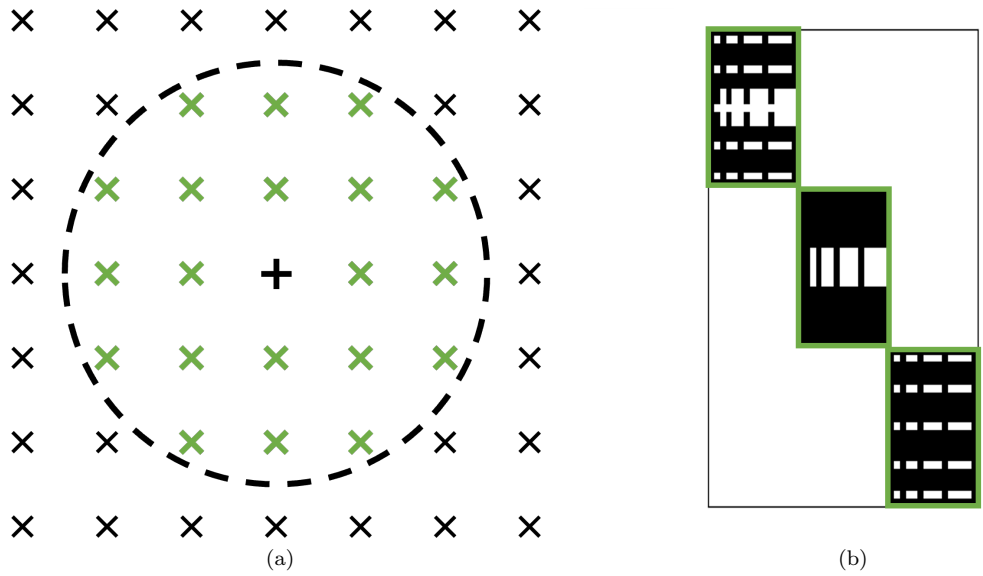
{A circular support region in free space and its corresponding matrix structure for the first-order acoustic wave equation (nonzero elements are black, zero elements are white). Staggered particle velocity subgrids are omitted to prevent excessive cluttering. Individual blocks are highlighted in green, and from top-left to bottom-right correspond to pressure, horizontal particle velocity, and vertical particle velocity. It is clear that this matrix can be split into three smaller systems which can be solved individually; each polynomial can be fitted independently of the others.
}

In free space, the stencil has an uninterrupted circular footprint comprising entirely of interior points. The corresponding matrix has a block-diagonal structure, as shown in Figure \ref{fig:matrix_no_boundary_annotated.png}, with each block corresponding to one of the fields. Each block can be inverted separately, meaning that the resultant polynomial extrapolations have no impact on one another. Note that splitting the matrix up in this manner considerably reduces the computational cost of finding the matrix inverse (or pseudoinverse as required), and it is thus desirable to do so where possible.

\multiplot*{2}{matrix_support_boundary.png,matrix_boundary_annotated.png}{width=0.8\columnwidth}
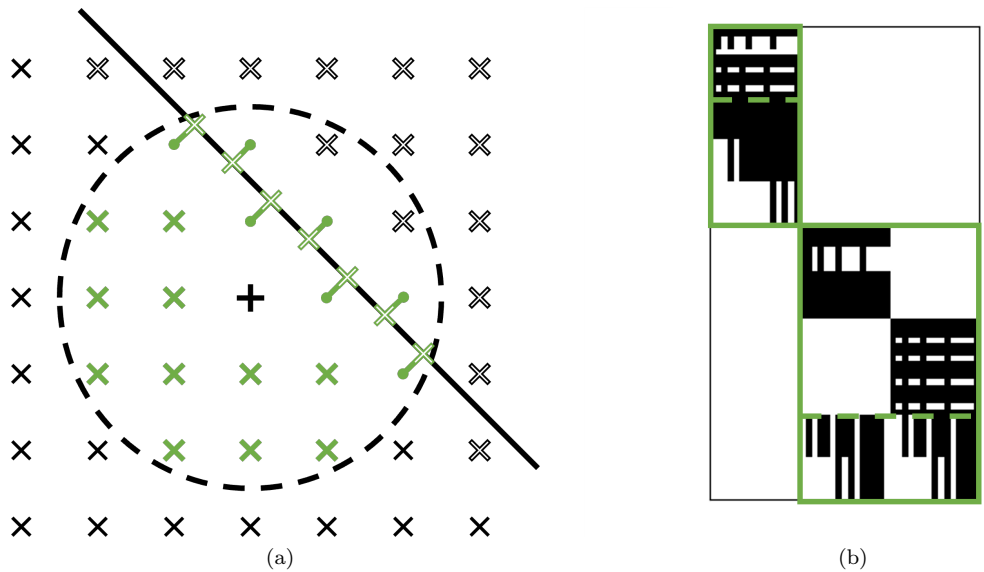
{A circular support region cut by a diagonal boundary and its corresponding matrix structure for the first-order acoustic wave equation (nonzero elements are black, zero elements are white). Staggered particle velocity subgrids are omitted to prevent excessive cluttering. Blocks are highlighted with solid green boxes, with the dashed green lines separating rows corresponding to interior points and those corresponding to boundary conditions. Note that the previously-separate particle velocity blocks have become merged, as the zero-divergence boundary condition (and its higher-order compatriots) span both particle velocity fields.
}

Inserting a 45\textdegree  free-surface boundary cutting across our support region, we see in Figure \ref{fig:matrix_boundary_annotated.png} that both the internal structure of the blocks and the overarching structure of the matrix are altered. Most notably, the inclusion of particle velocity boundary conditions has led to the merging of the corresponding blocks, meaning that they must be solved together, and resulting polynomial extrapolations will be dependent on both fields. As the pressure boundary conditions only concern the pressure field, this block retains its independence as before, although some rows are lost due to corresponding to now-external points, and several additional rows are added by the introduction of boundary conditions.

Whilst the matrices corresponding to such support regions will have consistently overdetermined blocks in free space, this may not be the case when boundaries are introduced. If a block becomes underdetermined, the radius of the support regions for functions contained therein can be expanded, thereby adding additional rows to the block until sufficient information is present.

\twocolumn

\bibliographystyle{seg}  
\bibliography{refs}

\end{document}